\documentclass[12pt]{article}

\usepackage{amsmath}
\allowdisplaybreaks
\usepackage{ifthen}\newboolean{draftmode}\setboolean{draftmode}{false}
\usepackage{exscale,relsize}
\usepackage{amsmath}
\usepackage{amsfonts}
\usepackage[colorlinks=true, linkcolor=blue]{hyperref}
\usepackage{amssymb}
\usepackage{calc}
\usepackage{theorem}
\usepackage{pifont}      
\usepackage{array}
\usepackage{color}
\usepackage{enumerate}
\usepackage{bbm}
\usepackage{graphicx}
\usepackage{subcaption}
\usepackage{caption}

\usepackage{algcompatible, algorithm}
\usepackage{mathtools}
\usepackage{tcolorbox}
\usepackage{soul}  

\ifthenelse{\boolean{draftmode}}{
    \usepackage[left]{showlabels}
    \usepackage[
        colorinlistoftodos,
        textsize=footnotesize, 
        color=yellow
    ]{todonotes} 
    
}{
    \usepackage[disable]{todonotes}  
}

\ifthenelse{\boolean{draftmode}}{
    \evensidemargin=\dimexpr\evensidemargin  + 6cm\relax 
    \oddsidemargin=\dimexpr\oddsidemargin + 6cm\relax
    \paperwidth=\dimexpr \paperwidth + 12cm\relax 
    \marginparwidth=\dimexpr \marginparwidth  + 6cm\relax
    \paperheight=\dimexpr \paperheight + 6cm\relax
}{
    
}

\newtheorem{theorem}{Theorem}[section]
\newtheorem{lemma}[theorem]{Lemma}

\newtheorem{proposition}[theorem]{Proposition}
\newtheorem{definition}[theorem]{Definition}

\theoremstyle{plain}{\theorembodyfont{\rmfamily}
\newtheorem{assumption}[theorem]{Assumption}}
\theoremstyle{plain}{\theorembodyfont{\rmfamily}
}
\theoremstyle{plain}{\theorembodyfont{\rmfamily}

\theoremstyle{plain}{\theorembodyfont{\rmfamily}
}
\theoremstyle{plain}{\theorembodyfont{\rmfamily}
\newtheorem{remark}[theorem]{Remark}}
\theoremstyle{plain}{\theorembodyfont{\rmfamily}
}

\numberwithin{equation}{section}

\newcommand{\RR}{\ensuremath{\mathbb R}}

\newcommand{\NN}{\ensuremath{\mathbb N}}

\newcommand{\dom}{\ensuremath{\operatorname{dom}}}

\newcommand{\pluss}{{\hskip1pt \raise1pt\vbox{\hrule width6pt \vskip1pt
\hrule width6pt}\kern-4pt{\lower1pt\hbox{\vrule height6pt \kern1pt\vrule
height6pt}}\hskip5pt}}

\newcommand{\argmin}{\mathop{\rm argmin}\limits}

\newcommand{\N}{\ensuremath{\mathbb N}}
\newcommand{\defeq}{\vcentcolon=}

\begin{document}

\title{{\fontfamily{ptm}\selectfont
        Relaxed Weak Accelerated Proximal Gradient Method: a Unified Framework for Nesterov's Accelerations
    }}

\author{
    Hongda Li
    \thanks{Department of Mathematics, I.K. Barber Faculty of Science,
    The University of British Columbia, Kelowna, BC Canada V1V 1V7. E-mail:  \texttt{alto@mail.ubc.ca}.}~ and ~Xianfu Wang
    \thanks{Department of Mathematics, I.K. Barber Faculty of Science,
    The University of British Columbia, Kelowna, BC Canada V1V 1V7. E-mail:  \texttt{shawn.wang@ubc.ca}.}
}

\date{\today}

\maketitle


\begin{abstract}
    \noindent
    This paper is devoted to the study of accelerated proximal gradient methods where the sequence that controls the momentum term doesn't follow Nesterov's rule.
    We propose a relaxed weak accelerated proximal gradient (R-WAPG) method, a generic algorithm that unifies the convergence results for strongly convex and convex problems where the extrapolation constant is characterized by a sequence that is much weaker than Nesterov's rule.
    Our R-WAPG provides a unified framework for several notable Euclidean variants of FISTA and verifies their convergences.
    In addition, we provide the convergence rate of strongly convex objective with a constant momentum term.
    Without using the idea of restarting, we also reformulate R-WAPG as ``Free R-WAPG" so that it doesn't require any parameter.
    Explorative numerical experiments were conducted to show its competitive advantages.
\end{abstract}


\section{Introduction}
    \subsection{Motivations}
        The Nesterov acceleration scheme, originally proposed in 1983 \cite{nesterov_method_1983}, is a celebrated first order method for solving the minimization problem.
        If the minimizer $x^*$ exists, then the optimality gap sequence $(F(x_k) - F(x^*))_{k\geq 1}$ decreases at a rate of $\mathcal O(1/k^2)$, a faster rate compared to gradient descent which is $\mathcal O(1/k)$.
        It's consider optimal in a certain sense, see Chapter 2 of Nesterov's book \cite{nesterov_lectures_2018} for details.
        However, things are not sunshine, rainbow and happily ever after because the above scheme is associated with the issue that if $F$ is $\mu > 0$ strongly convex, the optimality gap sequence $(F(x_k) - F(x^*))_{k\geq 1}$ oscillates and converges sub-linearly, making it worse than the gradient descent. See page 9 of Su et al. \cite{su_differential_2016}.
        \par
        This issue motivated a vast amount of literatures aims at improving, interpreting, and extending this method.
        Restarting is a popular solution to address the issue of obtaining better convergence rate when the objective function is strongly convex.
        Beck and Teboulle \cite{beck_fast_2009} mitigate the issue by restarting and showed that it has an $\mathcal O(1/k^2)$ convergence still, but it performs empirically better.
        See \cite{aujol_parameter-free_2024}, \cite{necoara_linear_2019} and references therein for recent advancements in restarting accelerated proximal gradient algorithms.
        \par
        In the convex setting, the variants that Nesterov's Accelerated Gradient are optimal to for cases when $\mu \ge 0$ and $\mu > 0$ are different.
        In this paper, we will show a unified framework for the convergence proof of both cases with a much relaxed condition for the recurrence of the momentum parameters used to extrapolate the iterates.
        For simplicity consider minimizing a differentiable convex function $F: \RR^n \rightarrow \RR$.
        Let $F$ be $L$-Lipschitz smooth and $\mu \ge 0$ strongly convex meaning that there exists $L \in \RR$ such that for all $x \in \RR^n, y\in \RR^n$ it has $\mu\Vert x - y\Vert\le \Vert \nabla F(x) - \nabla F(y)\Vert \le L \Vert x - y \Vert$.
        We will list the following example algorithm of Nesterov's Accelerated Gradient to refer to while introducing the context and prior works.
        Initialize $x_1 = y_1$ and $t_0 = 1$, the algorithm finds $(x_k)_{k \ge 1}$ for all $k \ge 1$ by:
        \begin{align}
            & x_{k + 1} := y_k - L^{-1}\nabla F(y_k),
            \\
            & t_{k + 1} := 1/2\left(1 + \sqrt{1 + 4t_{k}^2}\right),
            \\
            & \theta_{k + 1} := (t_{k} - 1)/t_{k + 1}, \label{eqn:example-algorithm}
            \\
            & y_{k + 1} := x_{k + 1} + \theta_{k + 1}(x_{k + 1} - x_k).
        \end{align}

        \subsection{Prior works}
            In \cite{apidopoulos_convergence_2018} Apidopoulos et al. showed that a relaxed choice of the sequence that doesn't strictly follow Nesterov's rule 
            still leads to the convergence of the momentum algorithm. In \cite{aujol_parameter-free_2024}, Aujol et al. gave a parameter free FISTA by adaptive restart and backtracking. 
           In Attouch, et al. \cite{attouch_first-order_2022}, Maulen and Peypouquet \cite{maulen_speed_2023}. the authors
            present the idea of Hessian Damping for stabilizing Nesterov's accelerated gradient method without/with restarting.

    \subsection{Contributions of this paper}
        Our are concerned with the composite model
        $\min_{x\in\RR^n}\{F(x)= f(x) + g(x)\},$ 
        where $g:\RR^n \rightarrow \overline\RR$ is a proper,  closed and convex function, and $f$ is an $L$-Lipschitz smooth and $\mu \ge 0$ strongly convex function.
        Our contributions are two folds: theoretical and practical.
        \par
        Inspired specifically by the technique of Nestrov's estimating sequence \cite{nesterov_lectures_2018}, firstly we present a unified framework of Accelerated Proximal Gradient (APG) which we call Relaxed Weak Accelerated Proximal Gradient (R-WAPG) in Section \ref{sec:rwapg-formulation-convergence}.
        It can upper bound $F(x_k) - F(x^*)$ for sequences $(t_k)_{k \ge 0}$ that doesn't follow Nesterov's update rule, and it unifies the convergence results of several Euclidean variants of FISTA.
        Secondly, as an alternative to restarting, we present the method of Free R-WAPG which estimates the modulus of convexity as the algorithm executes. Extensive numerical experiments illustrate that our method performs well empirically.
        \par
        \textbf{A summary of our main results follows. }
        Nesterov's acceleration extrapolates $y_{k + 1} := x_{k + 1} + \theta_{k + 1}(x_{k + 1} - x_k)$ where $\theta_{k + 1} := (t_{k} - 1)/t_{k + 1} \in (0, 1)$ is the ``momentum''.
        The choices for $\theta_k$ vary for different variants of the accelerated proximal gradient algorithm.
        In Chambolle and Dossal \cite{chambolle_convergence_2015}, it has $t_k := (k + a - 1)/a$ for all $a > 2$ which gives weak convergence of the iterates $x_k$ in Hilbert space.
        In Chapter 10 of Beck's Book \cite{beck_first-order_2017}, a variant called V-FISTA can achieve the faster linear convergence rate: $\mathcal O((1 - \sqrt{\mu/L})^k)$ on the optimality gap for $\mu > 0$ strongly convex $f$.
        V-FISTA has $\theta_k := (\sqrt{\kappa} - 1)/(\sqrt{\kappa} + 1)$ where $\kappa := L/\mu$.
        \par
        We relax the traditional choice of the sequence $(\theta_k)_{k\geq 1}$ in equation \eqref{eqn:example-algorithm} and give an upper bound of the optimal gap.
        Let $(\alpha_k)_{k \ge0}, (\rho_k)_{k \ge 0}$ be two sequences that satisfy
        \begin{align*}
            \alpha_0 &\in (0, 1],
            \\
            \alpha_k &\in (\mu/L, 1) \quad (\forall k \ge 1),
            \\
            \rho_k &:= \frac{\alpha_{k + 1}^2 - (\mu/L)\alpha_{k + 1}}{(1 - \alpha_{k + 1})\alpha_k^2} \quad \forall (k \ge 0).
        \end{align*}
        Our first main result shows that if $\theta_{k + 1} := (\rho_k\alpha_k(1 - \alpha_k)/(\rho_k\alpha_k^2 + \alpha_{k + 1}))$, using the R-WAPG we proposed in Definition \ref{def:wapg} with Proposition \ref{prop:wapg-convergence}, \ref{prop:r-wapg-momentum-repr}, we can show that the gap $F(x_k) - F(x^*)$ is bounded by:
        \begin{align*}
            \mathcal O\left(
                \left(
                    \prod_{i = 0}^{k - 1} \max(1, \rho_{i})
                \right)
                \prod_{i = 1}^{k} \left(1  - \alpha_i\right)
            \right).
        \end{align*}
        Our second main result shows that for any $a \ge 2$ the choice of sequence $\alpha_k := a/(a + k)$ results in $\rho_k \geq 1$ for all $k \in \N$ such that R-WAPG reduces to a variant of FISTA proposed in Chambolle and Dossal \cite{chambolle_convergence_2015}, and we are able to show the same convergence rate in Theorem \ref{thm:r-wapg-on-cham-doss}.
        When $\mu := 0, \rho_k := 1$, R-WAPG reduces perfectly to FISTA by Beck \cite{beck_first-order_2017}; when $\mu > 0, \rho_k := 1$, it reduces to the V-FISTA by Beck \cite{beck_first-order_2017}.
        In Theorem \ref{thm:fixed-momentum-fista}, we demonstrate that R-WAPG framework gives a linear convergence claim for all fixed momentum method where $\alpha_k := \alpha \in (\mu/L, 1)$ and  $F$ is $\mu > 0$ strongly convex.
        Finally, we present three equivalent forms of R-WAPG in Section \ref{sec:rwapg-equiv-repr} that are comparable to notable Euclidean variants of FISTA and beyond. This shows the enormous flexibility of our R-WAPG framework.
        \par
        Our practical contribution is an algorithm called ``Parameter Free R-WAPG'' (See Algorithm \ref{alg:free-rwapg}), and is is inspired by the proof of Proposition \ref{prop:wapg-convergence}.
        The algorithm is parameter free because it doesn't require knowing $L, \mu$ in advance, and it determines the value of $\theta_k$ by estimating the local concavity using iterates $y_{k}, y_{k + 1}$ from the Bregman divergence of $f$ with minimal computational cost.
        We conduct ample amount of numerical experiments to show that it has a favorable convergence rate in practice and behaves similarly to the FISTA with monotone restart.
        \par
        \textbf{The remainder of this paper is organized as follows. }
        In Section~\ref{sec:preliminaries}, we formulate the composite optimization
        problem considered in the paper and give a proximal inequality connecting the objective function,
        the proximal gradient
        operator and gradient mapping.
        Section \ref{sec:stepwise-stuff} proves 
        a key inequality based on a single generic iterative step that acts similar to a descent lemma; see Proposition \ref{prop:stepwise-lyapunov}. It is crucial for the convergence analysis of relaxed weak accelerated proximal gradient (R-WAPG) studied later.
        Section \ref{sec:rwapg-formulation-convergence} formulates the full R-WAPG and defines the R-WAPG sequence. Using Proposition \ref{prop:stepwise-lyapunov} from Section~\ref{sec:stepwise-stuff}
         and the R-WAPG sequence, 
        we provide the convergence rate of R-WAPG.
        Section \ref{sec:rwapg-equiv-repr} focuses on three equivalent representations of the R-WAPG algorithm that are comparable to instances of APG found in the literatures.
        Section \ref{sec:rwapg-literatures} formulates FISTA, and V-FISTA sequences as the instances of R-WAPG sequences.
        This section proves the convergence rate of several variants of FISTA using the equivalent forms introduced in Section \ref{sec:rwapg-equiv-repr} and the convergence rate developed in Section \ref{sec:rwapg-formulation-convergence}.
        In Section \ref{sec:free-rwapg} we give a parameter free version of R-WAPG algorithm and showcase its numerical experiments for regression, LASSO. We conclude with some discussions and perspectives of our results in Section~\ref{s:discussion}.

\section{Preliminaries}\label{sec:preliminaries}
    Throughout, we make the following assumption unless specified.
    \begin{assumption}[smooth and nonsmooth additives]\label{ass:smooth-nsmooth-additive}
        \;
        \begin{enumerate}
            \item Define $F := f + g$.
            \item $f: \RR^n \rightarrow \RR$ is $L$ Lipschitz smooth and $\mu \ge 0$ strongly convex.
            \item $g: \RR^n \rightarrow \overline \RR$ is proper, closed and convex. The extended real is defined as $\overline \RR := \RR \cup \{\infty\}$.
            \item Minimizer exists for the optimization problem: $F^* = \min_x \left\lbrace f(x) + g(x)\right\rbrace$.
        \end{enumerate}
    \end{assumption}
    For $x \in \RR^n, y \in \RR^n$, we define:
    \begin{align}
        \widetilde{\mathcal M}^{L^{-1}}
        (x; y)
        &:=
        g(x) + f(y) + \langle \nabla f(y), x - y\rangle
        + \frac{L}{2}\Vert x - y\Vert^2,  \label{eqn:pg-model-func}
        \\
        \mathcal M^{L^{-1}}(x; y) &:= F(x) + \frac{L}{2}\Vert x - y\Vert^2.
        \label{eqn:pp-model-func}
    \end{align}
    Define the proximal gradient operator and gradient mapping operator $T_L, \mathcal G_L$:
    \begin{align*}
        T_L(y) &:= \argmin_{x \in \RR^n} \left\lbrace
            g(x) + \langle \nabla f(y), x\rangle + L/2\Vert x - y\Vert^2
        \right\rbrace = [I + L^{-1}\partial g]^{-1}\circ [I - L^{-1}\nabla f](y),
        \\
        \mathcal G_L
        &:= L(I - T_L).
    \end{align*}
    $T_L$ is a single-valued mapping, and it has its domain on the entire $\RR^n$.
    Finally, define the Bregman divergence of $f$:
    \begin{align*}
        D_f(x, y): \RR^n \times \RR^n \rightarrow \RR
        \defeq (x, y)\mapsto f(x) - f(y) - \langle \nabla f(y), x - y\rangle.
    \end{align*}
    The following lemma provides important properties of $\mathcal M^{L^{-1}}(\cdot; y)$ and $ \widetilde{\mathcal M}^{L^{-1}}(\cdot; y)$.
    \begin{lemma}\label{lemma:pg-envelope}
        With $\mathcal M^{L^{-1}}, \widetilde{\mathcal M}^{L^{-1}}$ as given by \eqref{eqn:pg-model-func}, \eqref{eqn:pp-model-func},
        we have for all $x \in \RR^n$, $y \in \RR^n$:
        \begin{align*}
            \widetilde{\mathcal M}^{L^{-1}}(x; y)
            &=
            \mathcal M^{L^{-1}}(x; y)- D_f(x, y).
        \end{align*}
    \end{lemma}
    \begin{proof}
        The proof is direct using algebra and the definitions.
        For $x \in \dom g, y \in \RR^n$, using the definition of $\widetilde{\mathcal M}^{L^{-1}}(x; y)$ it has:
        \begin{align*}
            \widetilde{\mathcal M}^{L^{-1}}(x; y)
            &=
            g(x) + f(y) + \langle \nabla f(y), x - y\rangle + \frac{L}{2}\Vert x - y\Vert^2
            \\
            &=
            g(x) + f(x) - f(x) + f(y)
            + \langle \nabla f(y), x - y\rangle + \frac{L}{2}\Vert x - y\Vert^2
            \\
            &=
            F(x) - D_f(x, y) + \frac{L}{2}\Vert x - y\Vert^2
            \\
            &= \mathcal M^{L^{-1}}(x; y) - D_f(x, y).
        \end{align*}
        The equality is trivially true when $x \in \RR^n \setminus \dom g$.
    \end{proof}
    \begin{theorem}[proximal inequality]\label{thm:prox-grad-ineq}
        Let $F$ be given by Assumption \ref{ass:smooth-nsmooth-additive}.
        For $x,y\in \RR^n$, we have:
        \begin{align*}
            F(x) - F(T_Ly) - \langle L(y - T_Ly), x - y\rangle
            - \frac{L}{2}\Vert y - T_L y\Vert^2
            - \frac{\mu}{2}\Vert x - y\Vert^2
            &\ge
            0.
        \end{align*}
    \end{theorem}
    \begin{proof}
        For simplicity, we write $T$ for $T_L$.
        The proof is direct algebra, and it starts by deriving an intermediate inequality from the quadratic growth condition of the model function.
        Then it uses the characterization of strong convexity to simplify it to obtain the desired results.
        \par
        $\widetilde{\mathcal M}^{L^{-1}}(\cdot; y)$ is $L$ strongly convex by Assumption \ref{ass:smooth-nsmooth-additive} and its definition.
        We see that it admits quadratic growth over the minimizer $\bar y = T(y)$, therefore for all $x \in \RR^n$ it follows that:
        {\smaller
        \begin{align*}
            0 &\le
            \widetilde{\mathcal M}^{L^{-1}}(x; y) -
            \widetilde{\mathcal M}^{L^{-1}}(\bar y; y)
            -
            \frac{L}{2}\Vert x - \bar y\Vert^2
            \\
            &=
            \left(
                \mathcal M^{L^{-1}}(x; y) - D_f(x, y)
            \right) -
            \mathcal M^{L^{-1}}(\bar y; y)
            -
            \frac{L}{2}\Vert x - \bar y\Vert^2
            + D_f(\bar y; y)
            \quad
            \text{(By Lemma \ref{lemma:pg-envelope})}
            \\
            &=
            \left(
                \mathcal M^{L^{-1}}(x; y)
                -
                \mathcal M^{L^{-1}}(\bar y; y)
            \right)
            -
            D_f(x, y)
            - \frac{L}{2}\Vert x - \bar y\Vert^2
            + D_f(\bar y; y)
            \\
            &=
            \left(
                F(x) - F(\bar y)
                +
                \frac{L}{2}\Vert x - y\Vert^2 -
                \frac{L}{2}\Vert \bar y - y\Vert^2
            \right)
            - D_f(x, y)
            + D_f(\bar y; y)
            - \frac{L}{2}\Vert x - \bar y\Vert^2
            \\
            &=
            \left(
                F(x) - F(\bar y)
                +
                \frac{L}{2}
                \left(
                    \Vert x - \bar y + \bar y - y\Vert^2
                    -
                    \Vert y - \bar y\Vert^2
                \right)
            \right)
            - D_f(x, y)
            + D_f(\bar y; y)
            - \frac{L}{2}\Vert x - \bar y\Vert^2
            \\
            &=
            \left(
                F(x) - F(\bar y)
                +
                \frac{L}{2}
                \left(
                    \Vert x - \bar y\Vert^2 +
                    2\langle x - \bar y, \bar y - y\rangle
                \right)
            \right)
            - D_f(x, y)
            + D_f(\bar y; y)
            - \frac{L}{2}\Vert x - \bar y\Vert^2
            \\
            & =
            \left(
                F(x) - F(\bar y) + \frac{L}{2}\Vert x - \bar y\Vert^2
                - L\langle  x - \bar y, y - \bar y\rangle
            \right)
            - D_f(x, y)
            + D_f(\bar y; y)
            - \frac{L}{2}\Vert x - \bar y\Vert^2
            \\
            &=
            F(x) - F(\bar y)
            - \langle L(y - \bar y), x - \bar y\rangle
            - D_f(x, y)
            + D_f(\bar y, y).
        \end{align*}
        }
        From $\mu \ge 0$ strong convexity of $f$, it has $D_f(x, y) \ge \mu/2\Vert x - y\Vert^2$ which is true for all $x, y \in \RR^n$.
        Consequently, it has for all $x \in \RR^n, y \in \RR^n$:
        \begin{align*}
            0 &\le
            F(x) - F(\bar y)
            - \langle L(y - \bar y), x - \bar y\rangle
            - D_f(x, y)
            + D_f(\bar y, y)
            \\
            &\le
            F(x) - F(\bar y)
            - \langle L(y - \bar y), x - \bar y\rangle
            - \frac{\mu}{2}\Vert x - y\Vert^2
            + \frac{L}{2}\Vert y - \bar y\Vert^2
            \\
            &=
            F(x) - F(\bar y)
            - \langle L(y - \bar y), x - y  + y - \bar y\rangle
            - \frac{\mu}{2}\Vert x - y\Vert^2
            + \frac{L}{2}\Vert y - \bar y\Vert^2
            \\
            &=
            F(x) - F(\bar y) - \langle L(y - \bar y), x - y\rangle - \frac{\mu}{2}\Vert x - y \Vert^2
            - \frac{L}{2}\Vert y - \bar y\Vert^2.
        \end{align*}
        The proof is done.
    \end{proof}
    \begin{remark}
        This proof works for all $L$ that is larger than the smallest Lipschitz modulus of $f$ and $\mu$ less than the strong convexity modulus of $f$.
    \end{remark}

\section{Stepwise formulation of weak accelerated proximal gradient}\label{sec:stepwise-stuff}
    The goal of this section is to provide some preparatory results for the R-WAPG algorithm described in Definition \ref{def:wapg} of the next section.
    Definition \ref{def:stepwise-wapg} below describes one iteration of the R-WAPG algorithm.
    It defines the procedure of generating $x_{k + 1}, v_{k + 1}$ given any $x_k, v_k$.
    The main result of this section is Proposition \ref{prop:stepwise-lyapunov}, which states the key inequality \eqref{ineq:stepwise-lya-2}  of a decreasing quantity involving both $F(x_k)$ and $F(x_{k + 1})$ at each iteration.
    \par
    The assumption below is crucial to describing the convergence rate of the R-WAPG, and it requires the convexity assumption on $F$.
    \begin{assumption}\label{ass:regret-gap}
        Given $x_k, y_k, v_k$ where $k \in \mathbb \NN$, we define the following quantities:
        \begin{align}
            g_k &\defeq L(y_k - T_L(y_k)),
            \label{eqn:grad-map}
            \\
            l_F(x; y_k) &\defeq F(T_L(y_k)) + \langle g_k, x - y_k\rangle + \frac{1}{2L}\Vert g_k\Vert^2,
            \label{eqn:lower-linearization}
            \\
            \epsilon_{k} &\defeq F(x_k) - l_F(x_k; y_k).
            \label{eqn:regret}
        \end{align}
    \end{assumption}

    \begin{definition}[stepwise weak accelerated proximal gradient]\label{def:stepwise-wapg}\;\\
        Assume $0 \le \mu < L$.
        Fix any $k \in \NN$.
        For any $(v_k, x_k), \alpha_k \in (0, 1), \gamma_k > 0$, let $\hat \gamma_{k + 1}$, and vectors $y_k, v_{k + 1}, x_{k + 1}$ be given by:
        \begin{align}
            \hat \gamma_{k + 1} &:= (1 - \alpha_k)\gamma_k + \mu \alpha_k, \label{eqn:stepwise-wapg-eqn1}
            \\
            y_k &:=
            (\gamma_k + \alpha_k \mu)^{-1}(\alpha_k \gamma_k v_k + \hat\gamma_{k + 1} x_k), \label{eqn:stepwise-wapg-eqn2}
            \\
            g_k &:= \mathcal G_L (y_k), \label{eqn:stepwise-wapg-eqn3}
            \\
            v_{k + 1} &:= \hat\gamma^{-1}_{k + 1}
            (\gamma_k(1 - \alpha_k) v_k - \alpha_k g_k + \mu \alpha_k y_k), \label{eqn:stepwise-wapg-eqn4}
            \\
            x_{k + 1} &:= T_L(y_k). \label{eqn:stepwise-wapg-eqn5}
        \end{align}
    \end{definition}
    The following lemma will simplify the discussion later.
    \begin{lemma}\label{lemma:ineq-q3}
        Let $v_k, x_k, y_k, v_{k + 1}, x_{k + 1}$ and $\alpha_k, \hat \gamma_{k + 1}, \gamma_{k}$ be given by Definition \ref{def:stepwise-wapg}.
        Then for any $x^* \in \RR^n$ we have
        \begin{align}
            - \alpha_k(v_k - x^*) - \frac{\alpha_k^2 \mu}{\hat \gamma}(y_k - v_k) - (x_k - y_k)
            &= \alpha_k(x^* - x_k).
            \tag{Q3}\label{eqn:Q3}
        \end{align}
    \end{lemma}
    \begin{proof}
        We start by showing that \eqref{eqn:Q1} and \eqref{eqn:Q2} below 
        follow from \eqref{eqn:stepwise-wapg-eqn2}, and the desired result can be derived naturally using algebra.
        By Definition \ref{def:stepwise-wapg}:
        \begin{align*}
            y_k - v_k &=
            \frac{\hat \gamma_{k + 1}}{\alpha_k \gamma_k}(x_k - y_k),
            \tag{Q1}\label{eqn:Q1}
            \\
            y_k - x_k &=
            \frac{\alpha_k \gamma_k}{\gamma_k + \alpha_k \mu}(v_k - x_k).
            \tag{Q2}\label{eqn:Q2}
        \end{align*}
        It says that $x_k, y_k, v_k$ lie on the same line because \eqref{eqn:Q1} indicates $y_k - v_k$ parallels to $-(y_k - x_k)$ and both vector shares the same head which anchors at $y_k$.
        \par
        To see \eqref{eqn:Q1}, we show that it's equivalent to \eqref{eqn:stepwise-wapg-eqn2} by considering:
        \begin{align*}
            y_k - v_k &=
            \frac{\hat \gamma_{k + 1}}{\alpha_k \gamma_k}(x_k - y_k)
            \\
            \iff
            -(\alpha_k \gamma_k \hat \gamma^{-1}_{k + 1} + 1)y_k
            &=
            - \alpha_k \gamma_k \hat \gamma^{-1}_{k + 1}v_k - x_k
            \\
            \iff
            y_k &=
            \frac{
                \alpha_k \gamma_k \hat \gamma_{k + 1}^{-1}v_k + x_k
            }{1 + \alpha_k \gamma_k \hat \gamma_{k + 1}^{-1}}
            =
            \frac{\alpha_k \gamma_k v_k + \hat \gamma_{k + 1} x_k}{\gamma_k + \alpha_k \mu}.
        \end{align*}
        To see \eqref{eqn:Q2}, consider \eqref{eqn:stepwise-wapg-eqn2}:
        \begin{align*}
            y_k &= (\gamma_k + \alpha_k \mu)^{-1}(\alpha_k \gamma_k v_k + \hat\gamma_{k + 1} x_k)
            \\
            \iff
            y_k - x_k &=
            (\gamma_k + \alpha_k \mu)^{-1}
            (\alpha_k \gamma_k v_k - (\gamma_k + \alpha_k \mu)x_k + \hat \gamma_{k + 1} x_k)
            \\
            \iff
            (\gamma_k + \alpha_k \mu)(y_k - x_k)
            &=
            \alpha_k\gamma _kv_k +
            (\hat \gamma_k - \gamma_k - \alpha_k \mu) x_k
            & \text{by \eqref{eqn:stepwise-wapg-eqn1}}
            \\
            &= \alpha_k \gamma_k v_k - \alpha_k \gamma_k x_k
            = \alpha_k \gamma_k(v_k - x_k)
            \\
            \iff
            y_k - x_k &=
            \frac{\alpha_k \gamma_k}{\gamma_k + \alpha_k \mu}(v_k - x_k).
        \end{align*}
        \par
        Using \eqref{eqn:Q1}, \eqref{eqn:Q2}, we can show \eqref{eqn:Q3}.
        Indeed,
        \begin{align*}
            &
            - \alpha_k(v_k - x^*) - \frac{\alpha_k^2 \mu}{\hat \gamma_{k + 1}}(y_k - v_k) - (x_k - y_k)
            \\
            & \underset{\eqref{eqn:Q1}}{=}
            -\alpha_k(v_k - x^*) -
            \frac{\alpha_k^2\mu}{\hat \gamma_{k + 1}}
            \frac{\hat \gamma_{k + 1}}{\alpha_k \gamma_k}(x_k - y_k)
            - (x_k - y_k)
            \\
            &=
            -\alpha_k(v_k - x^*) -
            \frac{\alpha_k \mu}{\gamma_k}(x_k - y_k)
            - (x_k - y_k)
            \\
            &=
            -\alpha_k(v_k - x^*) -
            \left(
                1 + \frac{\alpha_k \mu}{\gamma_k}
            \right)(x_k - y_k)
            \\
            &\underset{\eqref{eqn:Q2}}{=}
            -\alpha_k(v_k - x^*) -
            \frac{\alpha_k \mu + \gamma_k}{\gamma_k}
            \frac{\alpha_k \gamma_k}{\gamma_k + \alpha_k \mu}(x_k - v_k)
            \\
            &=
            -\alpha_k(v_k - x^*)
            - \alpha_k(x_k - v_k)
            = \alpha_k(x^* - x_k).
        \end{align*}
    \end{proof}

    Our main result is this section comes as follows.
    \begin{proposition}[stepwise Lyapunov]\label{prop:stepwise-lyapunov}\;\\
        Let $k \in \mathbb Z_+$, $R_k \in \RR$.
        Given any $v_k, x_k$ and $\gamma_k > 0$ and $v_{k + 1}, x_{k + 1}, y_k, \hat \gamma_{k + 1}, \alpha_k$ that satisfy Definition \ref{def:stepwise-wapg}.
        With $\epsilon_k, y_k$ given by Assumption \ref{ass:regret-gap}, define:
        \begin{align}\label{eqn:stepwise-lya-1}
            R_{k + 1}
            \defeq
            \frac{1}{2}\left(
                L^{-1} - \frac{\alpha_k^2}{\hat \gamma_{k + 1}}
            \right)\Vert g_k\Vert^2
            +
            (1 - \alpha_k)
            \left(
                \epsilon_k + R_k +
                \frac{\mu\alpha_k\gamma_k}{2\hat \gamma_{k + 1}}
                \Vert v_k - y_k\Vert^2
            \right).
        \end{align}
        Then for all $x^* \in \RR^n$, we have:
        {\small
        \begin{align}\label{ineq:stepwise-lya-2}
            F(x_{k + 1}) - F(x^*) + R_{k + 1} + \frac{\hat \gamma_{k + 1}}{2}\Vert v_{k + 1} - x^*\Vert^2
            &\le
            (1 - \alpha_k)
            \left(
                F(x_k) - F(x^*) + R_k + \frac{\gamma_{k}}{2}\Vert v_k - x^*\Vert^2
            \right).
        \end{align}
        }
    \end{proposition}
    \begin{proof}
        For notation simplicity we drop the subscript $k$, $k + 1$ on $\gamma_{k}, \hat \gamma_{k + 1}$ because they are fixed throughout the proof.
        Proceeding stepwise, we will show the following intermediate results:
        \begin{enumerate}
            \item
            \begin{align}\tag{Step 1}\label{eqn:stepwise-lya-proof-step-1}
                \begin{split}
                    F(x_{k + 1}) + R_{k + 1}
                    &=
                    F(x_k) - \epsilon_k - \langle  g_k, x_k - y_k\rangle
                    - \frac{\alpha_k^2}{2\hat \gamma}\Vert g_k\Vert^2
                    \\
                    &\quad
                        +
                        (1 - \alpha_k)
                        \left(
                            \epsilon_k + R_k +
                            \frac{\mu\alpha_k\gamma}{2\hat \gamma}
                            \Vert v_k - y_k\Vert^2
                        \right).
                \end{split}
            \end{align}
            \item
            \begin{align}\tag{Step 2}\label{eqn:stepwise-lya-proof-step-2}
                \begin{split}
                    \frac{\hat \gamma}{2}\Vert v_{k + 1} - x^* \Vert^2
                    &=
                    \frac{(1 - \alpha_k)\gamma + \mu \alpha_k}{2} \Vert v_k - x^*\Vert^2
                    \\
                    & \quad
                        +
                        \left\langle g_k,
                            - \alpha_k(v_k - x^*)
                            - \frac{\alpha_k^2\mu}{\hat \gamma}(y_k - v_k)
                        \right\rangle
                    \\
                    & \quad
                        + \frac{\alpha_k^2}{2\hat \gamma}\Vert g_k\Vert^2
                        + \frac{\alpha_k^2 \mu^2}{2\hat \gamma}\Vert y_k - v_k\Vert^2
                        + \langle v_k - x^*, \mu\alpha_k(y_k - v_k)\rangle.
                \end{split}
            \end{align}
            \item Adding both side of Equations \eqref{eqn:stepwise-lya-proof-step-1}, \eqref{eqn:stepwise-lya-proof-step-2}, it follows:
            \begin{align}\tag{Step 3}\label{eqn:stepwise-lya-proof-step-3}
                \begin{split}
                    F(x_{k + 1}) + R_{k + 1} +
                    \frac{\hat \gamma}{2}\Vert v_{k + 1} - x^*\Vert^2
                    &=
                    F(x_k) - \alpha_k\left(
                        \epsilon_k + \langle g_k, x_k - x^*\rangle
                        - \frac{\mu}{2}\Vert y_k - x^*\Vert^2
                    \right)
                    \\
                    & \quad
                    +
                    (1 - \alpha_k)\left(
                        R_k + \frac{\gamma}{2}\Vert v_k - x^*\Vert^2
                    \right).
                \end{split}
            \end{align}
        \end{enumerate}
        We will prove the intermediate results above later.
        Building upon these intermediate results, we are able to claim \eqref{ineq:stepwise-lya-2}. Indeed,
         subtracting $F(x^*)$ on \eqref{eqn:stepwise-lya-proof-step-3} gives:
        \begin{align}
            \begin{split}\label{eqn:stepwise-lya-proof-eqn3.2}
                & F(x_{k + 1}) - F(x^*) + R_{k + 1} +
                \frac{\hat \gamma}{2}\Vert v_{k + 1} - x^*\Vert^2
                \\
                &=
                F(x_k) - F(x^*) - \alpha_k\left(
                    \epsilon_k + \langle g_k, x_k - x^*\rangle
                    - \frac{\mu}{2}\Vert y_k - x^*\Vert^2
                \right)
                +
                (1 - \alpha_k)\left(
                    R_k + \frac{\gamma}{2}\Vert v_k - x^*\Vert^2
                \right)
                \\
                &= (1 - \alpha_k)(F(x_k) - F(x^*))
                + \alpha_k\left(
                    F(x_k) - F(x^*) - \epsilon_k - \langle g_k, x_k - x^*\rangle + \frac{\mu}{2}\Vert y_k - x^*\Vert^2
                \right)
                \\ &\quad
                    +
                    (1 - \alpha_k)\left(
                        R_k + \frac{\gamma}{2}\Vert v_k - x^*\Vert^2
                    \right).
            \end{split}
        \end{align}
        Focusing on the second term, we simplify the multiplier inside:
        {\small
        \begin{align}
        \begin{split}\label{ineq:stepwise-lya-proof-ineq4}
            & F(x_k) - F(x^*) - \epsilon_k - \langle g_k, x_k - x^*\rangle + \frac{\mu}{2}\Vert y_k - x^*\Vert^2
            \\
            &=
            F(x_k) - F(x^*) - \left(
                F(x_k) - F(T_L y_k) - \langle g_k, x_k - y_k\rangle - \frac{1}{2L}\Vert g_k\Vert^2
            \right)- \langle g_k, x_k - x^*\rangle + \frac{\mu}{2}\Vert y_k - x^*\Vert^2
            \\
            &= F(T_L y_k) - F(x^*) + \langle g_k, x^* - y_k\rangle + \frac{\mu}{2}\Vert y_k - x^*\Vert^2
            + \frac{1}{2L}\Vert g_k\Vert^2 \underset{\text{Theorem }\ref{thm:prox-grad-ineq}}{\le} 0.
        \end{split}
        \end{align}
        }
        Then, by $\alpha_k > 0$, it has:
        {\small
        \begin{align*}
            F(x_{k + 1}) - F(x^*) + R_{k + 1} +
            \frac{\hat \gamma}{2}\Vert v_{k + 1} - x^*\Vert^2
            &\le
            (1 - \alpha_k)\left(
                F(x_k) - F(x^*) + R_k + \frac{\gamma}{2}\Vert v_k - x^*\Vert^2
            \right).
        \end{align*}
        } 
        \emph{Proof of \eqref{eqn:stepwise-lya-proof-step-1}}:
        Consider the first and the third term on LHS of \eqref{ineq:stepwise-lya-2} summed up:
        \begin{align}\label{eqn:stepwise-lya-proof-eqn1}
            F(x_{k + 1}) &\underset{\eqref{eqn:regret}}{=}
            F(x_k) - \epsilon_k - \langle  g_k, x_k - y_k\rangle - \frac{1}{2L}\Vert g_k\Vert^2,
            \notag
            \\
            R_{k + 1}
            &\underset{\eqref{eqn:stepwise-lya-1}}{=}
            \frac{1}{2}\left(
                L^{-1} - \frac{\alpha_k^2}{\hat \gamma}
            \right)\Vert g_k\Vert^2
            +
            (1 - \alpha_k)
            \left(
                \epsilon_k + R_k +
                \frac{\mu\alpha_k\gamma}{2\hat \gamma}
                \Vert v_k - y_k\Vert^2
            \right),
            \notag
            \\
            \implies
            F(x_{k + 1}) + R_{k + 1}
            &=
            F(x_k) - \epsilon_k - \langle  g_k, x_k - y_k\rangle
            - \frac{\alpha_k^2}{2\hat \gamma}\Vert g_k\Vert^2
            \notag
            \\
            &\quad
                +
                (1 - \alpha_k)
                \left(
                    \epsilon_k + R_k +
                    \frac{\mu\alpha_k\gamma}{2\hat \gamma}
                    \Vert v_k - y_k\Vert^2
                \right).
        \end{align}
        \emph{Proof of \eqref{eqn:stepwise-lya-proof-step-2}}:
        By definition of $v_{k + 1}$, it has:
        \begin{align}\label{eqn:stepwise-lya-proof-eqn2}\begin{split}
            \frac{\hat \gamma}{2}\Vert v_{k + 1} - x^* \Vert^2
            &=
            \frac{\hat \gamma}{2}\Vert
                \hat \gamma^{-1}
                (
                    \gamma(1 - \alpha_k)v_k -
                    \alpha_k g_k + \mu \alpha_k y_k
                )
                - x^*
            \Vert^2
            \\
            &=
            \frac{\hat \gamma}{2}
            \Vert
                \hat \gamma^{-1}
                (
                \hat \gamma v_k + \mu \alpha_k(y_k - v_k)
                    - \alpha_k g_k
                )
                - x^*
            \Vert^2
            \\
            &=
            \frac{\hat \gamma}{2}
            \Vert
                v_k + \hat \gamma^{-1} \mu \alpha_k (y_k - v_k)
                - \hat \gamma^{-1}\alpha_k g_k
                - x^*
            \Vert^2
            \\
            &=
            \frac{\hat \gamma}{2}
            \Vert v_k - x^*\Vert^2
            +
            \frac{\alpha_k^2}{2\hat \gamma}\Vert \mu(y_k - v_k) - g_k\Vert^2
            \\ &\quad
                +
                \langle v_k - x^*, \mu \alpha_k(y_k - v_k) - \alpha_k g_k\rangle
            \\
            &=
            \frac{(1 - \alpha_k)\gamma + \mu \alpha_k}{2} \Vert v_k - x^*\Vert^2
            \\ &\quad
                +
                \frac{\alpha_k^2}{2\hat \gamma}
                \Vert \mu(y_k - v_k) - g_k\Vert^2
                +
                \langle v_k - x^*, \mu \alpha_k(y_k - v_k) - \alpha_k g_k\rangle.
        \end{split}\end{align}
        The first equality above comes by substituting the definition of $v_{k + 1}$ given by \eqref{eqn:stepwise-wapg-eqn4}; the second equality simplifies using $\hat \gamma = (1 - \alpha_k)\gamma + \mu \alpha_k$ given by \eqref{eqn:stepwise-wapg-eqn1} and  because:
        \begin{align*}
            \gamma(1 - \alpha_k) v_k &=
            (\hat \gamma  - \mu \alpha_k)v_k
            = \hat \gamma v_k - \mu\alpha_k v_k
            \\
            \iff
            \gamma(1 - \alpha_k) v_k + \mu \alpha_k y_k
            &=
            \hat \gamma v_k + \mu \alpha_k(y_k - v_k).
        \end{align*}
        Focusing on the last two terms by the end of expression \eqref{eqn:stepwise-lya-proof-eqn2}, they can be written as:
        \begin{align}\label{eqn:stepwise-lya-proof-eqn2.1}\begin{split}
            \frac{\alpha^2_k}{2\hat \gamma}
            \Vert \mu(y_k - v_k) - g_k\Vert^2
            & =
            \frac{\alpha_k^2\mu}{\hat \gamma}
            \left(
                \frac{\mu}{2}\Vert y_k - v_k\Vert^2
                - \langle y_k - v_k, g_k\rangle
            \right)
            + \frac{\alpha_k^2}{2\hat \gamma}\Vert g_k\Vert^2,
            \\
            \langle v_k - x^*, \mu \alpha_k(y_k - v_k) - \alpha_k g_k\rangle
            &=
            \langle v_k - x^*, \mu\alpha_k(y_k - v_k)\rangle
            - \alpha_k \langle v_k - x^*, g_k\rangle.
        \end{split} 
        \end{align}
        Adding them gives:
        {\small
        \begin{align*}
            & \quad
            \frac{\alpha^2_k}{2\hat \gamma}
            \Vert \mu(y_k - v_k) - g_k\Vert^2
            +
            \langle v_k - x^*, \mu \alpha_k(y_k - v_k) - \alpha_k g_k\rangle
            \\
            &=
            \left\langle g_k,
                - \alpha_k(v_k - x^*)
                - \frac{\alpha_k^2\mu}{\hat \gamma}(y_k - v_k)
            \right\rangle
            + \frac{\alpha_k^2}{2\hat \gamma}\Vert g_k\Vert^2
            + \frac{\alpha_k^2 \mu^2}{2\hat \gamma}\Vert y_k - v_k\Vert^2
            + \langle v_k - x^*, \mu\alpha_k(y_k - v_k)\rangle.
        \end{align*}
        }
        With the above \eqref{eqn:stepwise-lya-proof-eqn2} simplifies to
        {\small
        \begin{align}\label{expr:stepwise-lya-expr2.2}
        \begin{split}
            \frac{\hat \gamma}{2}\Vert v_{k + 1} - x^*\Vert^2
            &=
            \frac{(1 - \alpha_k)\gamma + \mu \alpha_k}{2} \Vert v_k - x^*\Vert^2
            +
            \left\langle g_k,
                - \alpha_k(v_k - x^*)
                - \frac{\alpha_k^2\mu}{\hat \gamma}(y_k - v_k)
            \right\rangle
            \\
            & \quad
                + \frac{\alpha_k^2}{2\hat \gamma}\Vert g_k\Vert^2
                + \frac{\alpha_k^2 \mu^2}{2\hat \gamma}\Vert y_k - v_k\Vert^2
                + \langle v_k - x^*, \mu\alpha_k(y_k - v_k)\rangle.
        \end{split}
        \end{align}
        }
        \emph{Proof of \eqref{eqn:stepwise-lya-proof-step-3}}:
        Adding the RHS of \eqref{expr:stepwise-lya-expr2.2}, \eqref{eqn:stepwise-lya-proof-eqn1} gives:
        \begin{align}
        \begin{split}
            &
            F(x_k) - \epsilon_k - \langle  g_k, x_k - y_k\rangle
            - \frac{\alpha_k^2}{2\hat \gamma}\Vert g_k\Vert^2
            + (1 - \alpha_k)
            \left(
                \epsilon_k + R_k +
                \frac{\mu\alpha_k\gamma}{2\hat \gamma}
                \Vert v_k - y_k\Vert^2
            \right)
            \\
            &\quad
                +
                \frac{(1 - \alpha_k)\gamma + \mu \alpha_k}{2}
                \Vert v_k - x^*\Vert^2
                +
                \left\langle g_k,
                    - \alpha_k(v_k - x^*)
                    - \frac{\alpha_k^2\mu}{\hat \gamma}(y_k - v_k)
                \right\rangle
            \\
            & \quad
                + \frac{\alpha_k^2}{2\hat \gamma}\Vert g_k\Vert^2
                + \frac{\alpha_k^2 \mu^2}{2\hat \gamma}\Vert y_k - v_k\Vert^2
                + \langle v_k - x^*, \mu\alpha_k(y_k - v_k)\rangle
            \\
            &=
            F(x_k) - \epsilon_k
            + \left\langle
                g_k,
                - \alpha_k(v_k - x^*)
                - \frac{\alpha_k^2\mu}{\hat \gamma}(y_k - v_k)
                - (x_k - y_k)
            \right\rangle
            \\
            &\quad
                + (1 - \alpha_k)
                \left(
                    \epsilon_k + R_k +
                    \frac{\mu\alpha_k\gamma}{2\hat \gamma}
                    \Vert v_k - y_k\Vert^2
                \right)
                +
                \frac{(1 - \alpha_k)\gamma + \mu \alpha_k}{2} \Vert v_k - x^*\Vert^2
            \\
            & \quad
                + \frac{\alpha_k^2 \mu^2}{2\hat \gamma}\Vert y_k - v_k\Vert^2
                + \langle v_k - x^*, \mu\alpha_k(y_k - v_k)\rangle
            \\
            &\underset{\eqref{eqn:Q3}}{=}
            F(x_k) - \epsilon_k
            + \alpha_k\left\langle
                g_k,
                x^* - x_k
            \right\rangle
            \\
            &\quad
                + (1 - \alpha_k)
                \left(
                    \epsilon_k + R_k +
                    \frac{\mu\alpha_k\gamma}{2\hat \gamma}
                    \Vert v_k - y_k\Vert^2
                \right)
                +
                \frac{(1 - \alpha_k)\gamma + \mu \alpha_k}{2}
                \Vert v_k - x^*\Vert^2
            \\
            & \quad
                + \frac{\alpha_k^2 \mu^2}{2\hat \gamma}\Vert y_k - v_k\Vert^2
                + \langle v_k - x^*, \mu\alpha_k(y_k - v_k)\rangle
        \end{split}\nonumber
        \\
        \begin{split}
            &=
            F(x_k) - \alpha_k\epsilon_k + \alpha_k\langle g_k, x^* - x_k\rangle
            +
            (1 - \alpha_k)\left(
                R_k + \frac{\gamma}{2}\Vert v_k - x^*\Vert^2
            \right)
            \\&\quad
                + \frac{(1 - \alpha_k)\mu\alpha_k\gamma}{2\hat \gamma}\Vert v_k - y_k\Vert^2
                + \frac{\mu \alpha_k}{2}\Vert v_k - x^*\Vert^2
            \\&\quad
                + \frac{\alpha_k^2 \mu^2}{2\hat\gamma}\Vert y_k - v_k\Vert^2
                + \langle v_k - x^*, \mu\alpha_k(y_k - v_k)\rangle.
        \end{split}\label{eqn:stepwise-lya-proof-eqn3}
        \end{align}
        On the first equality, coefficients of $\Vert g_k\Vert^2$ cancels out to zero and the inner product containing $g_k$ are grouped
        using \eqref{eqn:Q3} from Lemma \ref{lemma:ineq-q3}.
        The last equality groups the coefficients of $\epsilon_k$ together.
        Continuing on \eqref{eqn:stepwise-lya-proof-eqn3} we have
        \begin{align}\label{eqn:stepwise-lya-proof-eqn3.1}
            &
            F(x_k) - \alpha_k(\epsilon_k + \langle g_k, x_k - x^*\rangle)
            +
            (1 - \alpha_k)\left(
                R_k + \frac{\gamma}{2}\Vert v_k - x^*\Vert^2
            \right)\nonumber
            \\&\quad
                + \frac{(1 - \alpha_k)\mu\alpha_k\gamma}{2\hat \gamma}\Vert v_k - y_k\Vert^2
                + \frac{\mu \alpha_k}{2}\Vert v_k - x^*\Vert^2 \nonumber
            \\&\quad
                + \frac{\alpha_k^2 \mu^2}{2\hat \gamma}\Vert y_k - v_k\Vert^2
                + \langle v_k - x^*, \mu\alpha_k(y_k - v_k)\rangle
            \nonumber\\
            &=
            F(x_k) - \alpha_k(\epsilon_k + \langle g_k, x_k - x^*\rangle)
            +
            (1 - \alpha_k)\left(
                R_k + \frac{\gamma}{2}\Vert v_k - x^*\Vert^2
            \right)
            \nonumber\\ &\quad
                +
                \left(
                    \frac{(1 - \alpha_k)\mu\alpha_k\gamma}{2\hat \gamma}
                    +
                    \frac{\alpha_k^2 \mu^2}{2\hat \gamma}
                \right)\Vert y_k - v_k\Vert^2
                + \frac{\mu \alpha_k}{2}\Vert v_k - x^*\Vert^2
                + \langle v_k - x^*, \mu\alpha_k(y_k - v_k)\rangle
            \nonumber\\
            & =
            F(x_k) - \alpha_k(\epsilon_k + \langle g_k, x_k - x^*\rangle)
            +
            (1 - \alpha_k)\left(
                R_k + \frac{\gamma}{2}\Vert v_k - x^*\Vert^2
            \right)
            \nonumber\\ &\quad
                +
                \frac{\mu \alpha_k}{2}\Vert y_k - v_k\Vert^2
                + \frac{\mu \alpha_k}{2}\Vert v_k - x^*\Vert^2
                + \langle v_k - x^*, \mu\alpha_k(y_k - v_k)\rangle
            \nonumber\\ &=
            F(x_k) - \alpha_k(\epsilon_k + \langle g_k, x_k - x^*\rangle)
            +
            (1 - \alpha_k)\left(
                R_k + \frac{\gamma}{2}\Vert v_k - x^*\Vert^2
            \right)
            \nonumber\\ &\quad
                +
                \frac{\mu\alpha_k}{2} \Vert y_k - x^*\Vert^2
            \nonumber
            \\&=
            F(x_k) - \alpha_k\left(
                \epsilon_k + \langle g_k, x_k - x^*\rangle
                - \frac{\mu}{2}\Vert y_k - x^*\Vert^2
            \right)
            +
            (1 - \alpha_k)\left(
                R_k + \frac{\gamma}{2}\Vert v_k - x^*\Vert^2
            \right).
        \end{align}
        On the first equality, coefficients of $\Vert y_k - v_k\Vert^2$ are grouped together, on the second to the third equality, its coefficient is simplified via:
        \begin{align*}
            \frac{(1 - \alpha_k)\mu\alpha_k\gamma}{2\hat \gamma} +
            \frac{\alpha_k^2 \mu^2}{2\hat \gamma}
            &=
            \frac{\mu\alpha_k}{2}\left(
                \frac{(1 - \alpha_k)\gamma_k + \alpha_k \mu}{\hat \gamma}
            \right)
            \\
            &\underset{\eqref{eqn:stepwise-wapg-eqn1}}{=} \frac{\mu\alpha_k}{2}\left(
                \frac{\hat \gamma}{\hat \gamma}
            \right) = \frac{\mu\alpha_k}{2}.
        \end{align*}
        We showed that \eqref{eqn:stepwise-lya-proof-eqn3.1} is the RHS of the sum \eqref{eqn:stepwise-lya-proof-step-1} + \eqref{eqn:stepwise-lya-proof-step-2}.
        Calling back to the LHS of the sum, it gives:
        \begin{align*}
            & F(x_{k + 1}) + R_{k + 1} +
            \frac{\hat \gamma}{2}\Vert v_{k + 1} - x^*\Vert^2
            \\
            &=
            F(x_k) - \alpha_k\left(
                \epsilon_k + \langle g_k, x_k - x^*\rangle
                - \frac{\mu}{2}\Vert y_k - x^*\Vert^2
            \right)
            +
            (1 - \alpha_k)\left(
                R_k + \frac{\gamma}{2}\Vert v_k - x^*\Vert^2
            \right).
        \end{align*}
        Hence \eqref{eqn:stepwise-lya-proof-step-3} is justified.
    \end{proof}
    \begin{remark}
        Given $y_k$, we can choose to increase $\mu_k = 2D_f(x^*, y_k)/\Vert y_k - x^*\Vert^2 \ge \mu$ which also works on \eqref{ineq:stepwise-lya-proof-ineq4}.
        This is true by the intermediate steps taken to prove Theorem \ref{thm:prox-grad-ineq}.
        $\mu$ is a pessimistic choice for the inequality above.
        But in general the choice of $\mu$ remains the strong convexity modulus or equivalently, any value that is smaller than the true strong convexity constant for claiming the convergence rate for all initial guesses.
    \end{remark}

\section{R-WAPG and its convergence rates}\label{sec:rwapg-formulation-convergence}
    In this section we propose \emph{Relaxed Weak Accelerated Proximal Gradient (R-WAPG)}, see Definition \ref{def:wapg}.
    The R-WAPG algorithm generates iterates $(x_k, y_k, v_k)$ and admits an upper bound on $F(x_k) - F^*$ described in Proposition \ref{prop:stepwise-lyapunov}.
    Definition \ref{def:rwapg-seq} introduces the R-WAPG sequences $(\alpha_k)_{k \ge 0}, (\rho_k)_{k \ge 0}$ which are crucial for our analysis.
    These sequences parameterize the R-WAPG algorithm, connect the step-wise formulation of R-WAPG (Definition \ref{def:stepwise-wapg}) and can guarantee the convergence of R-WAPG in Proposition \ref{prop:wapg-convergence}.
    In the next section, they continue to play an 
    essential role in describing several equivalent forms of the R-WAPG algorithm, and their corresponding convergence claims.
    \begin{definition}[R-WAPG sequences]\label{def:rwapg-seq}\;\\
        Assume $0 \le \mu < L$.
        The sequences $(\alpha_k)_{k \ge 0}, (\rho_k)_{k \ge 0}$ are valid for R-WAPG if all the following holds:
        \begin{align*}
            \alpha_0 &\in (0, 1],
            \\
            \alpha_k &\in (\mu/L, 1) \quad (\forall k \ge 1),
            \\
            \rho_k &:= \frac{\alpha_{k + 1}^2 - (\mu/L)\alpha_{k + 1}}{(1 - \alpha_{k + 1})\alpha_k^2} \quad \forall (k \ge 0).
        \end{align*}
        We call $(\alpha_k)_{k \ge 0}, (\rho_k)_{k \ge 0}$ as the {R-WAPG sequences}.
    \end{definition}
    Not all properties of the R-WAPG sequences are used in the convergence proof, some other properties such as $\alpha_k \in (\mu/L, 1)$ plays a role much later on in Section \ref{sec:rwapg-equiv-repr}.
    The following lemma states two consequences of the R-WAPG sequences, which will play a crucial role in the convergence proof immediately after.
    \begin{lemma}\label{lemma:r-wapg-seq-consequences}
        Let sequence $(\alpha_k)_{k\ge 0}$ and $(\rho_k)_{k \ge0}$ be given by Definition \ref{def:rwapg-seq}. Then 
        \begin{enumerate}
            \item\label{i:rho} $\rho_k > 0$ for all $k \ge 0$.
            \item\label{i:alpha} $\forall k \ge 1: L\alpha_k^2 = (1 - \alpha_k)L\rho_{k - 1}\alpha_{k - 1}^2 +\mu\alpha_k$.
        \end{enumerate}
    \end{lemma}
    \begin{proof}
        \ref{i:rho} is obvious by $\alpha_k \in (\mu/L, 1)$.
        \ref{i:alpha} is direct by definition of $\rho_k$ which has for all $k \ge 1$:
        \begin{align*}
            (1 - \alpha_k)L \rho_{k - 1}\alpha_{k - 1}^2 + \mu\alpha_k
            &= (1 - \alpha_k)L \left(
                \frac{\alpha_{k}^2 - (\mu/L)\alpha_{k}}{(1 - \alpha_{k})\alpha_{k - 1}^2}
            \right)\alpha_{k - 1}^2 + \mu\alpha_k
            \\
            &= (L\alpha_k^2 - \mu \alpha_k) + \mu \alpha_k = L\alpha_k^2.
        \end{align*}
    \end{proof}
    \par
    We now give Relaxed Weak Accelerated Proximal Gradient in details.
    \begin{definition}[relaxed weak accelerated proximal gradient (R-WAPG)]\label{def:wapg}\;\\
        Choose any $x_1 \in \RR^n, v_1 \in \RR^n$.
        Let $(\alpha_k)_{k \ge0}, (\rho_k)_{k \ge 0}$ be given by Definition \ref{def:rwapg-seq}.
        The algorithm generates a sequence of vector $(y_k, x_{k + 1}, v_{k + 1})_{k \ge 1}$ for $k\ge 1$ by the procedures:
        \begin{tcolorbox}
            For $k=1, 2, 3, \ldots$
            \begin{align*}
                \gamma_k &\defeq \rho_{k -1}L\alpha_{k - 1}^2,
                \\
                \hat \gamma_{k + 1} & \defeq (1 - \alpha_k)\gamma_k + \mu \alpha_k = L\alpha_k^2,
                \\
                y_k &\defeq
                (\gamma_k + \alpha_k \mu)^{-1}(\alpha_k \gamma_k v_k + \hat\gamma_{k + 1} x_k),
                \\
                g_k &\defeq \mathcal G_L (y_k),
                \\
                v_{k + 1} &\defeq
                \hat\gamma^{-1}_{k + 1}
                (\gamma_k(1 - \alpha_k) v_k - \alpha_k g_k + \mu \alpha_k y_k),
                \\
                x_{k + 1} &\defeq T_L y_k.
            \end{align*}
        \end{tcolorbox}
    \end{definition}
    \begin{remark}\label{remark:rwapg-def}
        For all $k \ge 1$, the algorithm chains together the sequence of $(\hat \gamma_{k+1})_{k \ge 1}$ and $(\gamma_k)_{k \ge 1}$ in Definition \ref{prop:stepwise-lyapunov} through the equalities:
        \begin{align*}
            \hat \gamma_{k + 1}
            &= L\alpha_k^2
            = (1 - \alpha_k)\gamma_k + \alpha_k \mu
            \\
            &= (1 - \alpha_k)L\rho_{k - 1}\alpha_{k - 1}^2 + \alpha_k \mu
            = (1 - \alpha_k)\rho_{k - 1}\hat \gamma_{k} + \alpha_k \mu.
        \end{align*}
        Therefore, the sequences $(\alpha_k)_{k \ge 0}, (\rho_k)_{k \ge 0}$ fits the ecurrence relation from Lemma \ref{lemma:r-wapg-seq-consequences}\ref{i:alpha}.
        \par
        Observe that if $\rho_k = 1$ for all $k\ge 0$, then the above algorithm is similar to (2.2.7) in Nesterov's book \cite{nesterov_lectures_2018} because $L\alpha_k^{2} = (1 - \alpha_k)L\alpha_{k - 1}^2 + \mu \alpha_k$.
    \end{remark}
    Here is the main result of this section.
    \begin{proposition}[R-WAPG convergence]\label{prop:wapg-convergence}\; \\
        Fix any arbitrary $x^* \in \RR^n, N \in \mathbb N$.
        Let vector sequence $(y_k, v_{k}, x_{k})_{k \ge 1}$ and R-WAPG sequences $\alpha_k, \rho_k$ be given by Definition \ref{def:wapg}.
        Let $\epsilon_k, y_k$ be given by Assumption \ref{ass:regret-gap}.
        Define $R_1 = 0$ and suppose that for $k = 1, 2, \ldots, N$, we have $R_k$ recursively given by:
        \begin{align*}
            R_{k + 1}
            :=
            \frac{1}{2}\left(
                L^{-1} - \frac{\alpha_k^2}{\hat \gamma_{k + 1}}
            \right)\Vert g_k\Vert^2
            +
            (1 - \alpha_k)
            \left(
                \epsilon_k + R_k +
                \frac{\mu\alpha_k\gamma_k}{2\hat \gamma_{k + 1}}
                \Vert v_k - y_k\Vert^2
            \right).
        \end{align*}
        Then for all $k = 1, 2, \ldots, N$:
        \begin{align*}
            & F(x_{k + 1}) - F(x^*) + \frac{L \alpha_k^2}{2}\Vert v_{k + 1} - x^*\Vert^2
            \\
            &\le
            \left(
                \prod_{i = 0}^{k - 1} \max(1, \rho_{i})
            \right)
            \left(
                \prod_{i = 1}^{k} \left(1  - \alpha_i\right)
            \right)
            \left(
                F(x_1) - F(x^*) + \frac{L\alpha_0^2}{2}\Vert v_1 - x^*\Vert^2
            \right).
        \end{align*}
    \end{proposition}

    \begin{proof}
        We prove by induction.
        Let's first consider $k = 1$, $\hat \gamma_2 = (1 - \alpha_1)\gamma_1 + \mu \alpha_1$, $R_1 = 0$.
        From the above definition, $R_2$ is the same as \eqref{eqn:stepwise-lya-1}.
        Invoking Proposition \ref{prop:stepwise-lyapunov}, we have:
        \begin{align*}
            & F(x_{2}) - F(x^*) + R_2 + \frac{L \alpha_1^2}{2}\Vert v_{2} - x^*\Vert^2
            \\
            & = F(x_{2}) - F(x^*) + R_2 + \frac{\hat \gamma_{2}}{2}\Vert v_{2} - x^*\Vert^2
            & (\hat \gamma_2 = L \alpha_1^2 , \text{ Definition \ref{def:wapg}})
            \\
            &\le
            (1 - \alpha_1)\left(
                F(x_1) - F(x^*) + R_1 + \frac{\gamma_1}{2}\Vert v_1 - x^*\Vert^2
            \right)
            & \quad \text{(By Proposition \ref{prop:stepwise-lyapunov})}
            \\
            &=
            (1 - \alpha_1)\left(
                F(x_1) - F(x^*) + R_1 + \frac{L\rho_{0}\alpha_{0}^2}{2}\Vert v_1 - x^*\Vert^2
            \right)
            \\
            &= \max\left(\rho_0, 1\right)
            (1 - \alpha_1)\left(
                F(x_1) - F(x^*) + R_1 + \frac{L\alpha_0^2}{2}\Vert v_1 - x^*\Vert^2
            \right).
        \end{align*}
        We had demonstrated the base case.
        Next, for all $k = 2, 3, \ldots, N$, from Remark \ref{remark:rwapg-def}:
        \begin{align*}
            \hat \gamma_{k + 1} = L\alpha_{k}^2
            &=(1 - \alpha_k)\rho_{k - 1}L\alpha_{k - 1}^2 + \mu\alpha_k
            \\
            &= (1 - \alpha_k)\gamma_k + \mu\alpha_k.
        \end{align*}
        So $\gamma_k, \hat \gamma_{k + 1}, \alpha_k$ fits Definition \ref{def:stepwise-wapg}.
        For all $k \ge 1$, $R_k$ satisfies \eqref{eqn:stepwise-lya-1}, hence unroll recursively using Proposition \ref{prop:stepwise-lyapunov}:
        {\small
        \begin{align*}
            &
            F(x_{k + 1}) - F^* + R_{k + 1} + \frac{\hat\gamma_{k + 1}}{2}\Vert v_{k + 1} - x^*\Vert^2
            \\
            &=
            F(x_{k + 1}) - F^* + R_{k + 1} + \frac{L \alpha_k^2}{2}\Vert v_{k + 1} - x^*\Vert^2
            \\
            &\le
            (1 - \alpha_k)
            \left(
                F(x_k) - F^* + R_k + \frac{\rho_{k - 1}L \alpha_{k - 1}^2}{2}\Vert v_k - x^*\Vert^2
            \right)
            \\
            &\le
            (1 - \alpha_k)
            \left(
                F(x_k) - F^* + R_k + \max(1, \rho_{k - 1})\frac{L \alpha_{k - 1}^2}{2}\Vert v_k - x^*\Vert^2
            \right)
            \\
            &\le
            \max(1, \rho_{k - 1})(1 - \alpha_k)
            \left(
                F(x_k) - F^* + R_k + \frac{L \alpha_{k - 1}^2}{2}\Vert v_k - x^*\Vert^2
            \right)
            \\
            &\le
            \left(
                \prod_{i = 0}^{k - 1} \max(1, \rho_{i})
            \right)
            \left(
                \prod_{i = 1}^{k} \left(1  - \alpha_i\right)
            \right)
            \left(
                F(x_1) - F^* + R_1 + \frac{L\alpha_0^2}{2}\Vert v_1 - x^*\Vert^2
            \right).
        \end{align*}
        }
        We had demonstrated in the inductive case for $k=1, 2, \ldots, N$.
        Additionally, for all $k = 1, 2, \ldots, N$ it has $R_{k + 1} \ge 0$ because:
        \begin{align*}
            R_{k + 1}
            &=
            \frac{1}{2}\left(
                L^{-1} - \frac{\alpha_k^2}{\hat \gamma_{k + 1}}
            \right)\Vert g_k\Vert^2
            +
            (1 - \alpha_k)
            \left(
                \epsilon_k + R_k +
                \frac{\mu\alpha_k\gamma_k}{2\hat \gamma_{k + 1}}
                \Vert v_k - y_k\Vert^2
            \right)
            \\
            &= (1 - \alpha_k)
            \left(
                \epsilon_k + R_k
                + \frac{\mu\alpha_k\gamma_k}{2\hat \gamma_{k + 1}}
                \Vert v_k - y_k\Vert^2
            \right)
            \\
            &\ge
            (1 - \alpha_k) R_k
            \\
            &\ge R_1 \prod_{i = 1}^{k} \left(1 - \alpha_i\right) = 0.
        \end{align*}
        Going from the left to the right on the first equality, we used the fact that $\hat \gamma_{k + 1} = L \alpha_{k}^2$.
        This makes coefficient of $\Vert g_k\Vert^2$ zero.
        The first inequality is by $\epsilon_k \ge 0$ and the non-negativity of the one remaining term.
        The last equality is by the assumption that $R_1 = 0$.
        Therefore:
        {\small
        \begin{align*}
            &
            F(x_{k + 1}) - F^* +
            \frac{L\alpha_k^2}{2}\Vert v_{k + 1} - x^*\Vert^2
            \\
            &\le
            \left(
                \prod_{i = 0}^{k - 1} \max(1, \rho_{i})
            \right)
            \left(
                \prod_{i = 1}^{k} \left(1  - \alpha_i\right)
            \right)
            \left(
                F(x_1) - F^* + \frac{L\alpha_0^2}{2}\Vert v_1 - x^*\Vert^2
            \right).
        \end{align*}
        }
    \end{proof}
    \begin{remark}
        The choice of $\rho_k = 1$ reproduces traditional choice of sequences for the Accelerated Proximal Gradient (APG) algorithm.
    \end{remark}

\section{Equivalent representations of R-WAPG}\label{sec:rwapg-equiv-repr}
    The goal of this section is to reduce Definition \ref{def:wapg} into various forms that are equivalent to what commonly appear in the literature.
    Variants of Accelerated Proximal Gradient algorithm such as FISTA, V-FISTA have different representations.
    Sometimes the same algorithm include two iterates $(x_k, y_k)$; sometimes they include three iterates $(x_k, y_k, z_k)$.
    In this section we simplify R-WAPG to these representations commonly found in the literature and showed that they are equivalent.
    These equivalent representations are listed in Definitions \ref{def:r-wapg-intermediate}, \ref{def:r-wapg-st-form} and \ref{def:r-wapg-momentum-form}.
    These forms are equivalent, and they are crucial for proving the convergence of existing variants found in the literature; see Section \ref{sec:rwapg-literatures}.
    \par
    Proposition \ref{prop:wapg-first-equivalent-repr} simplifies Definition \ref{def:wapg} and finds a representation without using auxiliary sequence $\gamma_k, \hat \gamma_k$.
    Definition \ref{def:r-wapg-intermediate} states the first simplified form of the R-WAPG algorithm which we call: ``R-WAPG intermediate form''.
    Following a similar pattern, Proposition \ref{prop:wagp-st-form}, \ref{prop:r-wapg-momentum-repr} demonstrates two more equivalent representations of the R-WAPG intermediate form (Definition \ref{def:r-wapg-intermediate}) which are formulated into Definition \ref{def:r-wapg-st-form}, \ref{def:r-wapg-momentum-form}.
    Convergence results from Proposition \ref{prop:wapg-convergence} applies to all these equivalent forms of R-WAPG.
    In brief, different equivalent reformulations are summarized as following:
    \begin{align*}
        &\text{Definition \ref{def:wapg}}   \iff
        \text{Definition \ref{def:r-wapg-intermediate}}  & \text{(By Proposition \ref{prop:wapg-first-equivalent-repr}).}
        \\
        & \text{Definition \ref{def:r-wapg-intermediate}}
        \iff \text{Definition \ref{def:r-wapg-st-form}} & \text{(By Proposition \ref{def:r-wapg-st-form}).}
        \\
        &
        \text{Definition \ref{def:r-wapg-st-form}}\implies
        \text{Definition }\ref{def:r-wapg-momentum-form}.
        & \text{(By Proposition \ref{prop:r-wapg-momentum-repr})}.
    \end{align*}
    We start with the following result on ``abstract similar triangle form'' to make later proofs cleaner.
    \begin{proposition}[abstract similar triangle form]\label{prop:abs-st-form}\;\\
        Given any initial $(x_1, v_1)$, and a sequence $(\tau_k)_{k \ge 1}, (\xi_k)_{k \ge 1}$ such that for all $k \ge 1$ $\tau_k \in (0, 1), \xi_k \in (0, 1)$, and the iterates $(y_k, z_{k + 1}, x_{k + 1})_{k \ge 1}$ satisfy recursively:
        \begin{align*}
            y_k &= (1 + \tau_k)^{-1}(v_k + \tau_k x_k),
            \\
            v_{k + 1} &= (1 + \xi_k)^{-1}(v_k + \xi_k y_k) - (1 + \xi_k)^{-1}\delta_k g_k,
            \\
            x_{k + 1} &= y_k - L^{-1} g_k.
        \end{align*}
        If $1 + \xi_k + \tau_k = L\delta_k\; \forall k \ge 1$, then for all $k \ge 1$ we have
        \begin{align}\label{eqn:abs-st-key}
            v_{k + 1} - x_{k + 1} = (1 + \xi_k)^{-1}\tau_k(x_{k + 1} - x_k),
        \end{align}
        which makes the algorithm a similar triangle form.
    \end{proposition}
    \begin{proof}
        We are interested in identifying the conditions required for the sequence of $\xi_k, \tau_k, \delta_k$ such that there exists $\theta_k$ satisfying:
        \begin{align*}
            v_{k + 1} - x_{k + 1}
            &= \theta_k(x_{k + 1} - x_k).
        \end{align*}
        To verify, we consider:
        \begin{align*}
            v_{k + 1} &=
            (1 + \xi_k)^{-1}(v_k + \xi_t y_k - \delta_k \mathcal G_L(y_k))
            \\
            &=
            (1 + \xi_k)^{-1}((1 + \tau_k)y_k - \tau_t x_k + \xi_k y_k - \delta_k \mathcal G_L(y_k))
            \\
            &=
            (1 + \xi_k)^{-1}((1 + \tau_k + \xi_k)y_k - \tau_k x_k - \delta_k \mathcal G_L(y_k))
            \\
            \iff
            v_{k + 1} - x_{k + 1}
            &=
            (1 + \xi_k)^{-1}((1 + \tau_k + \xi_k)y_t - \tau_k x_k - \delta_t \mathcal G_L(y_k))
            - (y_k - L^{-1}\mathcal G_L(y_k))
            \\
            &=
            (1 + \xi_k)^{-1}(\tau_ky_k - \tau_k x_k - \delta_k \mathcal G_L(y_k))
            + L^{-1}\mathcal G_L(y_k)
            \\
            &=
            (1 + \xi_k)^{-1}
            \left(
                \tau_ty_k - \tau_t x_k + (L^{-1}(1 + \xi_k) - \delta_k) \mathcal G_L(y_k)
            \right)
            \\
            &=
            (1 + \xi_k)^{-1}\tau_k
            \left(
                y_k - x_k +
                \tau_k^{-1}(L^{-1}(1 + \xi_k) - \delta_k) \mathcal G_L(y_k)
            \right).
        \end{align*}
        Going between the first and second inequality we used $v_k = (1 + \tau_k)y_k - \tau_k x_k$ which is rearranged $y_k = (1 + \tau_k)^{-1}(v_k + \tau_k x_k)$.
        The RHS is can be verified through
        \begin{align*}
            x_{k + 1} - x_k &=
            y_t - L^{-1}\mathcal G_L(y_k) - x_k
            \\
            &= (y_k - x_k) - L^{-1}\mathcal G_L(y_k).
        \end{align*}
        It necessitates the condition:
        \begin{align*}
            \tau_k^{-1}(L^{-1}(1 + \xi_k) - \delta_k)
            &= - L^{-1}
            \\
            (1 + \xi_k) - L\delta_k
            &=
            - \tau_k
            \\
            1 + \xi_k + \tau_k
            &=
            L\delta_k.
        \end{align*}
        Thus,
        \begin{align*}
            v_{k + 1} - x_{k + 1} &=
            (1 + \xi_k)^{-1}\tau_t
            \left(y_k - x_k - L^{-1}\mathcal G_L(y_k)\right)
            =
            (1 + \xi_k)^{-1}\tau_k(x_{k + 1} - x_k).
        \end{align*}
    \end{proof}
    \subsection{Equivalent representations of R-WAPG}
        This section lists three equivalent representations of the R-WAPG algorithm.
        They are comparable to existing APG algorithms in the literature such as Exercise 12.1 in Ryu, Yin \cite{ryu_large-scale_2022}, similar triangle form in Lee et al. \cite{lee_geometric_2021} and Ahn Sra \cite{ahn_understanding_2022}, and momentum form of (2.2.19) in Nesterov \cite{nesterov_lectures_2018}.
        Specific instances of Accelerated Proximal Gradient algorithm that has the same form as the Definition \ref{def:r-wapg-intermediate}, Definition \ref{def:r-wapg-st-form} and Definition \ref{def:r-wapg-momentum-form} in the literature are stated in the remarks that follow the definitions.
        \begin{definition}[R-WAPG intermediate form]\label{def:r-wapg-intermediate}\;\\
            Assume $\mu < L$ and let $(\alpha_k)_{k \ge 0}, (\rho_k)_{k \ge 0}$ given by Definition \ref{def:rwapg-seq}.
            Initialize any $x_1, v_1$ in $\RR^n$.
            For $k \ge 1$, the algorithm generates sequence of vector iterates $(y_{k}, v_{k + 1}, x_{k + 1})_{k \ge 1}$ by the procedures:
            \begin{tcolorbox}
                For $k = 1, 2, \ldots$
                \begin{align*}
                    & y_{k} \defeq
                    \left(
                        1 + \frac{L - L\alpha_{k}}{L\alpha_{k} - \mu}
                    \right)^{-1}
                    \left(
                        v_{k} +
                        \left(\frac{L - L\alpha_{k}}{L\alpha_{k} - \mu} \right) x_{k}
                    \right),
                    \\
                    & x_{k + 1} \defeq
                    y_k - L^{-1} \mathcal G_L (y_k),
                    \\
                    & v_{k + 1} \defeq
                    \left(
                        1 + \frac{\mu}{L \alpha_k - \mu}
                    \right)^{-1}
                    \left(
                        v_k +
                        \left(\frac{\mu}{L \alpha_k - \mu}\right) y_k
                    \right) - \frac{1}{L\alpha_{k}}\mathcal G_L (y_k).
                \end{align*}
            \end{tcolorbox}
        \end{definition}
        \begin{remark}
            This form of APG is rarely identified in the literatures.
            The closest algorithm that fits the form but with $\mu = 0$ is Chapter 12 of in Ryu and Yin's Book \cite{ryu_large-scale_2022}, right after Theorem 17.
            We created this form which makes the math that follows simpler.
            The inspiration of using this as an intermediate representation was inspired by solving Exercise 12.1 in the same Ryu and Yin's Book.
        \end{remark}
        \begin{definition}[R-WAPG similar triangle form]\label{def:r-wapg-st-form} \; \\
            Given any $(x_1, v_1)$ in $\RR^n$.
            Assume $\mu < L$.
            Let the sequence $(\alpha_k)_{k \ge 0}, (\rho_k)_{k\ge 0}$ be given by Definition \ref{def:rwapg-seq}.
            For $k \ge 1$, the algorithm generates sequences of vector iterates $(y_k, v_{k + 1}, x_{k + 1})_{k \ge 1}$ by the procedures:
            \begin{tcolorbox}
                For $k=1, 2, \ldots $
                \begin{align*}
                    & y_k \defeq
                    \left(
                        1 + \frac{L - L\alpha_k}{L\alpha_k - \mu}
                    \right)^{-1}
                    \left(
                        v_k +
                        \left(\frac{L - L\alpha_k}{L\alpha_k - \mu} \right) x_k
                    \right),
                    \\
                    & x_{k + 1} \defeq
                    y_k - L^{-1} \mathcal G_L (y_k),
                    \\
                    & v_{k + 1} \defeq
                    x_{k + 1} + (\alpha_k^{-1} -1)(x_{k + 1} - x_k).
                \end{align*}
            \end{tcolorbox}
        \end{definition}
        \begin{remark}
            The word similar triangle form can be traced back to several literatures.
            The term ``Method of Similar Triangle" was used for Algorithm (6.1.19) in Nesterov's book \cite{nesterov_lectures_2018}, but without the necessary graphical illustrations to clarify it.
            Equation (2), (3), (4) in \cite{chambolle_convergence_2015} is a similar triangle formulation of FISTA with $\mu = 0$.
            To see graphical visualization on why such term is used to describe the APG algorithm in the literatures, see
            Sections (3.1), (4.1) in Lee et al. \cite{lee_geometric_2021} and Ahn and Sra \cite{ahn_understanding_2022}.
        \end{remark}
        \begin{definition}[R-WAPG momentum form]\label{def:r-wapg-momentum-form}
            Given any $y_1 = x_1 \in \RR^n$, and sequences $(\rho_k)_{k \ge 0}, (\alpha_k)_{k\ge 0}$ Definition \ref{def:rwapg-seq}.
            The algorithm generates iterates $x_{k + 1}, y_{k + 1}$ For $k = 1, 2, \cdots $ by the procedures:
            \begin{tcolorbox}
                For $k=1, 2,\ldots $
                \begin{align*}
                    & x_{k + 1} \defeq y_k - L^{-1}\mathcal G_L (y_k),
                    \\
                    &
                    y_{k + 1} \defeq
                    x_{k + 1} +
                    \frac{\rho_k\alpha_k(1 - \alpha_k)}{\rho_k\alpha_k^2 + \alpha_{k + 1}}(x_{k + 1} - x_k).
                \end{align*}
            \end{tcolorbox}
            In the special case where $\mu = 0$, the momentum term can be represented without relaxation parameter $\rho_k$:
            $$
                (\forall k \ge 1)\quad \frac{\rho_k\alpha_k(1 - \alpha_k)}{\rho_k\alpha_k^2 + \alpha_{k + 1}}
                = \alpha_{k + 1}(\alpha_k^{-1} - 1).
            $$
        \end{definition}
        \begin{remark}
            This format fits with (2.2.19) in Nesterov's book \cite{nesterov_lectures_2018}, however, the sequence $(\alpha_k)_{k \ge 0}$ would be given by a different rule.
            See Theorem \ref{thm:r-wapg-on-cham-doss} and Lemma \ref{lemma:inverted-fista-seq} for a specific choice of $(\alpha_k)_{k \ge0}, (\rho_k)_{ k\ge 0}$ such that this equivalent form of R-WAPG is in fact two possible variants of the FISTA algorithm.
        \end{remark}

    \subsection{Derivations of equivalent representations of R-WAPG}
        \begin{proposition}[R-WAPG intermediate form]\label{prop:wapg-first-equivalent-repr}\;\\
            If the sequence $(y_k, v_k, x_k)_{k \ge 1}$ is produced by R-WAPG (Definition \ref{def:wapg}),
            then the iterates can be expressed without $(\gamma_k)_{k \ge1},(\hat \gamma_k)_{k \ge 2}$,  namely, for all $k\ge 1$,
            \begin{align}
                & y_{k} =
                \left(
                    1 + \frac{L - L\alpha_{k}}{L\alpha_{k} - \mu}
                \right)^{-1}
                \left(
                    v_{k} +
                    \left(\frac{L - L\alpha_{k}}{L\alpha_{k} - \mu} \right) x_{k}
                \right), \label{eqn:rwapg-first-equiv-form-eqn-1}
                \\
                & x_{k + 1} =
                y_k - L^{-1} \mathcal G_L (y_k),
                \\
                & v_{k + 1} =
                \left(
                    1 + \frac{\mu}{L \alpha_k - \mu}
                \right)^{-1}
                \left(
                    v_k +
                    \left(\frac{\mu}{L \alpha_k - \mu}\right) y_k
                \right) - \frac{1}{L\alpha_{k}}\mathcal G_L (y_k).
                \label{eqn:rwapg-first-equiv-form-eqn-2}
            \end{align}
        \end{proposition}
        \begin{proof}
            {We will show \eqref{eqn:rwapg-first-equiv-form-eqn-1} first.}
            For all $k \ge 1$, by R-WAPG (Definition \ref{def:wapg}), it has:
            \begin{align*}
                y_{k} &=
                (\gamma_k + \alpha_k \mu)^{-1}
                (\alpha_k \gamma_k v_k + \hat \gamma_{k + 1}x_k)
                \\
                &=
                (\hat \gamma_{k + 1} + \alpha_k \gamma_k)^{-1}
                (\alpha_k \gamma_k v_k + \hat \gamma_{k + 1}x_k)
                \\
                &=
                \left(
                    \frac{\hat \gamma_{k + 1}}{\alpha_k\gamma_k} + 1
                \right)^{-1}
                \left(
                    v_k + \frac{\hat \gamma_{k + 1}}{\alpha_k \gamma_k} x_k
                \right)
                =
                \left(
                    \frac{L\alpha_k^2}{\alpha_k\gamma_k} + 1
                \right)^{-1}
                \left(
                    v_k + \frac{L\alpha_k^2}{\alpha_k \gamma_k} x_k
                \right)
                \\
                &=
                \left(
                    \frac{L\alpha_k}{\gamma_k} + 1
                \right)^{-1}
                \left(
                    v_k + \frac{L\alpha_k}{ \gamma_k} x_k
                \right)
                \\
                &=
                \left(
                    1 + \frac{L - L \alpha_k}{L \alpha_k - \mu}
                \right)^{-1}
                \left(
                    v_k +
                    \frac{L - L \alpha_k}{L \alpha_k - \mu} x_k
                \right).
            \end{align*}
            From the left to right of the second equality, we used $\hat \gamma_{k + 1} = (1 - \alpha_k)\gamma_k + \alpha_k\mu$ which comes from \eqref{eqn:stepwise-wapg-eqn2}.
            Going from the left to the right of the second last equality, we did the following:
            \begin{align*}
                L\alpha_k^2 &=
                (1 - \alpha_k)\gamma_k + \alpha_k \mu
                \\
                \iff
                L \alpha_k^2 - \alpha_k\mu &=
                (1 - \alpha_k)\gamma_k
                \\
                \iff
                \gamma_k/L
                &=
                \frac{L \alpha_k^2 - \alpha_k\mu}{L (1 - \alpha_k)}
                \\
                \iff
                L/\gamma_k
                &=
                \frac{L (1 - \alpha_k)}{L \alpha_k^2 - \alpha_k\mu}
                \\
                \iff
                L\alpha_k/\gamma_k
                &=
                \frac{L - L\alpha_k}{L\alpha_k - \mu}.
            \end{align*}
            On the third $\iff$, we can assume $\alpha_k \neq \mu/L\;  \forall k \ge 1$ because from Definition \ref{def:rwapg-seq}: $\alpha_k \in (\mu/L, 1)$ for all $k \ge 1$.
            \par
            {Next, we show \eqref{eqn:rwapg-first-equiv-form-eqn-2}.}
            From the Definition \ref{def:wapg}, for all $k \ge 1$, $v_{k + 1}$ it follows that:
            \begin{align*}
                v_{k + 1} &=
                \hat \gamma_{k + 1}^{-1}
                ((1 - \alpha_k)\gamma_k v_k + \mu\alpha_k y_k)
                - \alpha_k\hat \gamma_{k + 1}^{-1}\mathcal G_L (y_k)
                \\
                &=
                ((1 - \alpha_k)\gamma_k + \alpha_k \mu)^{-1}
                \left(
                    (1 - \alpha_k)\gamma_k v_k + \mu\alpha_k y_k
                \right)
                - \alpha_k\hat \gamma_{k + 1}^{-1}\mathcal G_L (y_k)
                \\
                &=
                \left(
                    1 + \frac{\alpha_k\mu}{(1 - \alpha_k)\gamma_k}
                \right)^{-1}
                \left(
                    v_k +
                    \frac{\alpha_k\mu}{(1 - \alpha_k)\gamma_k} y_k
                \right)
                - \alpha_k\hat \gamma_{k + 1}^{-1}\mathcal G_L (y_k)
                \\
                &=
                \left(
                    1 + \frac{\alpha_k \mu}{L \alpha_k^2 - \alpha_k \mu}
                \right)^{-1}
                \left(
                    v_k +
                    \frac{\alpha_k \mu}{L \alpha_k^2 - \alpha_k \mu} y_k
                \right)
                - \frac{1}{L\alpha_{k}}\mathcal G_L (y_k).
            \end{align*}
            Going from the left to the right of the second equality, we substitute $\hat \gamma_{k + 1} = (1 - \alpha_k)\gamma_k + \mu\alpha_k$.
            At the end, recall that for all $k \ge 1$, it has $\hat \gamma_{k + 1} = L \alpha_k^2 = (1 - \alpha_k)\gamma_k + \alpha_k \mu$, so:
            \begin{align*}
                (1 - \alpha_k)\gamma_k
                &=
                \hat \gamma_{k + 1} - \mu \alpha_k
                =
                L\alpha_{k}^2 - \alpha_k\mu.
            \end{align*}
            The proof is now complete.
            This form doesn't have $\rho_k, \gamma_k, \hat \gamma_k$ in it.
        \end{proof}

        \begin{proposition}[R-WAPG similar triangle form]\label{prop:wagp-st-form}\;\\
            Let iterates $(y_k, x_{k}, v_{k})_{k \ge 1}$ and sequence $(\alpha_k, \rho_k)_{k \ge 0}$ be given by Definition \ref{def:r-wapg-intermediate}.
            Then for all $k \ge 1$, iterates $y_k, x_{k + 1}, v_{k + 1}$
            satisfy:
            \begin{align}
                y_{k} &=
                \left(
                    1 + \frac{L - L\alpha_{k}}{L\alpha_{k} - \mu}
                \right)^{-1}
                \left(
                    v_{k} +
                    \left(\frac{L - L\alpha_{k}}{L\alpha_{k} - \mu} \right) x_{k}
                \right),
                \label{eqn:rwapg-st-form-eqn-1}
                \\
                x_{k + 1} &=
                y_k - L^{-1} \mathcal G_L (y_k),
                \\
                v_{k + 1} &= x_{k + 1} + (\alpha_k^{-1} - 1)(x_{k + 1} - x_k).
                \label{eqn:rwapg-st-form-eqn-3}
            \end{align}
        \end{proposition}
        \begin{proof}
            From Definition \ref{def:r-wapg-intermediate}, define $(\tau_k, \xi_k, \delta_k)_{k \ge 1}$ in Proposition \ref{prop:abs-st-form} to be
            \begin{align}
                (\forall k \ge 1) \quad
                \tau_k &= \frac{L(1 - \alpha_k)}{L\alpha_k - \mu},
                \label{eqn:rwapg-st-form-proof-tau}
                \\
                (\forall k \ge 1)\quad
                \xi_k &= \frac{\mu}{L \alpha_k - \mu},
                \\
                (\forall k \ge 1)\quad
                \delta_k &\defeq \frac{1 + \xi_k}{L\alpha_k}.
            \end{align}
            Then, it fits into the format of abstract similar triangle form Proposition \ref{prop:abs-st-form}, and it follows that $v_{k + 1} - x_{k + 1} = (1 + \xi_k)^{-1}(x_{k + 1} - x_k)\; \forall k \ge 1$.
            We shall verify the hypothesis that $L\delta_k = 1 + \tau_k + \xi_k$ later.
            Substituting $\xi_k$ into \eqref{eqn:abs-st-key} gives:
            \begin{align*}
                v_{k + 1} &=
                x_{k + 1} + \left(
                    1 + \frac{\mu}{L\alpha_k - \mu}
                \right)^{-1}\left(
                    \frac{L(1 - \alpha_k)}{L\alpha_k - \mu}
                \right)(x_{k + 1} - x_k)
                \\
                &=
                x_{k + 1} + \left(
                    \frac{L\alpha_k}{L\alpha_k - \mu}
                \right)^{-1}\left(
                    \frac{L(1 - \alpha_k)}{L\alpha_k - \mu}
                \right)(x_{k + 1} - x_k)
                \\
                &=
                x_{k + 1} + \left(
                    \frac{L\alpha_k - \mu}{L\alpha_k}
                \right)\left(
                    \frac{L - L\alpha_k}{L\alpha_k - \mu}
                \right)(x_{k + 1} - x_k)
                \\
                &= x_{k + 1} + (\alpha_k^{-1} - 1)(x_{k + 1} - x_k).
            \end{align*}
            It's easy to see that substituting $\tau_k$ as given by \eqref{eqn:rwapg-st-form-proof-tau} into $y_k = (1 + \tau_k)^{-1}(v_k + \tau_k x_k)$ as in Proposition \ref{prop:abs-st-form} produces \eqref{eqn:rwapg-st-form-eqn-1}.
            It remains to show that $L\delta_k = 1 + \tau_k + \xi_k$.
            To see this, we have for all $k\ge 1$:
            \begin{align*}
                1 + \tau_k + \xi_k &=
                1 + \frac{L(1 - \alpha_k)}{L \alpha_k - \mu}
                + \frac{\mu}{L \alpha_k - \mu}
                \\
                &=
                1 + \frac{L - L \alpha_k + \mu}{L\alpha_k - \mu}
                =
                \frac{L - L \alpha_k + \mu + L \alpha_k - \mu}{L\alpha_k - \mu}
                \\
                &= \frac{L}{L\alpha_k - \mu}.
            \end{align*}
            Next, it also has for all $k \ge 1$:
            \begin{align*}
                L\delta_k = \frac{1 + \xi_k}{\alpha_k}
                &=
                \frac{1 + \frac{\mu}{L\alpha_k - \mu}}{\alpha_k}
                =
                \frac{\frac{L\alpha_k - \mu + \mu}{L \alpha_k - \mu}}{\alpha_k}
                =
                \frac{L}{L\alpha_k - \mu}.
            \end{align*}

        \end{proof}

        \begin{proposition}[R-WAPG momentum form]\label{prop:r-wapg-momentum-repr}
            \;\\
            Let sequence $(\alpha_k, \rho_k)_{k \ge 0}$ and iterates $(x_k, v_k, y_k)_{k\ge 1}$ be given by the R-WAPG intermediate form (Definition \ref{def:r-wapg-st-form}).
            Then it has:
            \begin{enumerate}
                \item For all $k \ge 1$, the iterates $x_{k + 1}, y_{k + 1}$ satisfy:
                \begin{align*}
                    x_{k + 1} &= y_k - L^{-1}\mathcal G_L (y_k),
                    \\
                    y_{k + 1} &=
                    x_{k + 1} +
                    \frac{\rho_k\alpha_k(1 - \alpha_k)}
                    {\rho_k\alpha_k^2 + \alpha_{k + 1}}(x_{k + 1} - x_k).
                \end{align*}
                \item $v_1 = x_1$ if and only if $y = x_1$.
                \item If in addition $\mu = 0$, then the momentum term in (i) admits a simpler representation:
                \begin{align*}
                    (\forall k \ge 1) \quad
                    \frac{\rho_k\alpha_k(1 - \alpha_k)}{\rho_k\alpha_k^2 + \alpha_{k + 1}}
                    & = \alpha_{k + 1}(\alpha_k^{-1} - 1).
                \end{align*}
            \end{enumerate}
        \end{proposition}
        \begin{proof}
            Start by considering the update rule for $v_k$ from Definition \ref{def:r-wapg-st-form} which has for all $k \ge 1$:
            \begin{align}
                v_{k + 1} &=
                x_{k + 1} + (\alpha_k^{-1} - 1)(x_{k + 1} - x_k)
                \notag
                \\
                \iff
                (L \alpha_{k + 1} - \mu)v_{k + 1}
                &=
                (L \alpha_{k + 1} - \mu)x_{k + 1} + (L\alpha_{k + 1} - \mu)(\alpha_k^{-1} - 1)(x_{k + 1} - x_k).
                \label{eqn:third-eqv-repr-proof-step1}
            \end{align}
            Simplifying $y_k$ firstly for $k = 1$, then for all $k \ge 2$.
            By Definition \ref{def:r-wapg-st-form}, it follows that for all $k \ge 1$:
            \begin{align}\label{eqn:third-eqv-repr-proof-step2}
                \begin{split}
                    y_k &=
                    \left(
                        1 + \frac{L - L\alpha_k}{L\alpha_k - \mu}
                    \right)^{-1}
                    \left(
                        v_k +
                        \left(\frac{L - L\alpha_k}{L\alpha_k - \mu} \right) x_k
                    \right)
                    \\
                    &=
                    \left(
                    \frac{L - \mu}{L\alpha_k - \mu}
                    \right)^{-1}
                    \left(
                        v_k +
                        \left(\frac{L - L\alpha_k}{L\alpha_k - \mu} \right) x_k
                    \right)
                    \\
                    &=
                    \frac{L\alpha_k - \mu}{L - \mu} v_k
                    +
                    \frac{L - L \alpha_k}{L - \mu} x_k
                    \\
                    &= (L - \mu)^{-1}((L \alpha_k - \mu)v_k + (L - L \alpha_k)x_k).
                \end{split}
            \end{align}
            For all $k\ge 1$, substitute RHS of \eqref{eqn:third-eqv-repr-proof-step1} into the above to obtain:
            {\small
            \begin{align*}
                y_{k + 1} &=
                (L - \mu)^{-1}((L\alpha_{k + 1} - \mu)v_{k + 1} + (L - L \alpha_{k + 1})x_{k + 1})
                \\
                &= (L - \mu)^{-1}
                \left(
                    (L\alpha_{k + 1} - \mu)x_{k + 1} +
                    (L\alpha_{k + 1} - \mu)(\alpha_k^{-1} - 1)(x_{k + 1} - x_k)
                    + (L - L \alpha_{k + 1})x_{k + 1}
                \right)
                \\
                &=
                (L - \mu)^{-1}
                \left(
                    (L - \mu)x_{k + 1} + (L\alpha_{k + 1} - \mu)(\alpha_k^{-1} - 1)(x_{k + 1} - x_k)
                \right)
                \\
                &= x_{k + 1} + \frac{(L\alpha_{k + 1} - \mu)(\alpha_k^{-1} - 1)}{L - \mu}(x_{k + 1} - x_k).
            \end{align*}
            }
            To show (ii), we show $y_1 = x_1$ when $v_1 = x_1$.
            Put $k = 1$ in \eqref{eqn:third-eqv-repr-proof-step1} to obtain:
            \begin{align*}
                y_1 &= (L - \mu)^{-1}((L\alpha_1 - \mu)x_1 + (L - L\alpha_1)x_1)
                \\
                &= (L - \mu)^{-1}((L\alpha_1 - \mu)x_1 + (L - L \alpha_1)x_1)
                \\
                &= (L - \mu)^{-1}((L - \mu)x_1) = x_1.
            \end{align*}
            Observe that, if $\mu = 0$, we can prove (iii) because:
            \begin{align*}
                \frac{(L\alpha_{k + 1} - \mu)(\alpha_k^{-1} - 1)}{L - \mu} =
                \frac{L\alpha_{k +1}(\alpha_k^{-1} - 1)}{L} = \alpha_{k +1}(\alpha_k^{-1} - 1).
            \end{align*}
            Therefore, the special case had been justified.
            To prove (ii), it remains to show that for all $k\ge 1$ the coefficient of $x_{k + 1} - x_k$ can be expressed without $\mu, L$:
            \begin{align*}
                \frac{(L\alpha_{k + 1} - \mu)(\alpha_k^{-1} - 1)}{L - \mu}
                &= \frac{\rho_k\alpha_k(1 - \alpha_k)}{\rho_k\alpha_k^2 + \alpha_{k + 1}}.
            \end{align*}
            Consider:
            \begin{align*}
                \frac{(L\alpha_{k + 1} - \mu)(\alpha_k^{-1} - 1)}{L - \mu}
                &= \frac{(L\alpha_{k + 1} - \mu)\alpha_k(1 - \alpha_k)}{\alpha_k^2(L - \mu)}
                \\
                &=
                \alpha_k(1 - \alpha_k)
                \left(
                    \frac{\alpha_k^2(L - \mu)}{L\alpha_{k + 1} - \mu}
                \right)^{-1}
                \\
                &= \alpha_k(1 - \alpha_k)
                \left(
                    \frac{L\alpha_k^2 - \mu\alpha_k^2}{L\alpha_{k + 1} - \mu}
                \right)^{-1}
                \\
                &=
                \alpha_k(1 - \alpha_k)
                \rho_k\left(
                    \frac{L\rho_k\alpha_k^2 - \mu\rho_k\alpha_k^2}{L\alpha_{k + 1} - \mu}
                \right)^{-1}
                \\
                &=
                \rho_k\alpha_k(1 - \alpha_k)
                \left(
                    \frac{(L\alpha_{k + 1} - \mu)(\rho_k\alpha_k^2 + \alpha_{k + 1})}
                    {L\alpha_{k + 1} - \mu}
                \right)^{-1}
                \\
                &= \frac{\rho_k\alpha_k(1 - \alpha_k)}{\rho_k\alpha_k^2 + \alpha_{k + 1}}.
            \end{align*}
            Going from the left to right on the fourth equality, we used the property of R-WAPG sequence $L\alpha_{k + 1}^2 = (1 - \alpha_{k + 1})\rho_kL\alpha_k^2 + \mu \alpha_{k + 1}$ to establish:
            \begin{align*}
                L \rho_k \alpha_k^2 - \mu \rho_k \alpha_k^2
                &=
                (1 - \alpha_{k + 1})L \rho_k \alpha_k^2 + \alpha_{k + 1} L \rho_k \alpha_k^2 - \mu \rho_k \alpha_k^2
                \\
                &=
                ((1 - \alpha_{k + 1})L \rho_k \alpha_k^2 + \mu \alpha_{k + 1}) - \mu\alpha_{k + 1} + \alpha_{k + 1} L \rho_k \alpha_k^2 - \mu \rho_k \alpha_k^2
                \\
                &= L \alpha_{k + 1}^2 - \mu\alpha_{k + 1} + \alpha_{k + 1}L\rho_k\alpha_k^2 - \mu \rho_k \alpha_k^2
                \\
                &=
                L\alpha_{k + 1}(\alpha_{k + 1} + \rho_k \alpha_k^2) - \alpha_{k + 1}\mu - \mu \rho_k \alpha_k^2
                \\
                &= (L \alpha_{k + 1} - \mu)(\alpha_{k + 1} + \rho_k \alpha_k^2).
            \end{align*}
        \end{proof}

\section{R-WAPG unifies existing acceleration schemes}\label{sec:rwapg-literatures}
    In addition to various equivalent forms of the R-WAPG algorithm, the R-WAPG sequences are much more flexible.
    They significantly generalize many existing sequences used in accelerated proximal gradient schemes.
    \par
    Using results from Section \ref{sec:rwapg-equiv-repr}, this section will demonstrate that several Euclidean variants of FISTA in the literatures reduce to our R-WAPG for specific R-WAPG sequences with additional constraints.
    As a result, Proposition \ref{prop:wapg-convergence} applies. By substituting specialized R-WAPG sequences $(\alpha_k)_{k \ge 0}$, we derive convergence rates consistent with results in the literature for these variants of FISTA (fast iterative shrinkage-thresholding algorithm).
    \par
    The table below shows that our R-WAPG convergence results unify convergence claims for various settings in the literature.
    The first row shows the generic convergence results of our R-WAPG.
    The rows after are algorithms found in the literature.
    First columns states the assumption made for parameter $\mu$, the second and third columns show the specific examples of $(\alpha_k)_{k \ge 0}, (\rho_k)_{k \ge0}$ specialized from Definition \ref{def:rwapg-seq}.
    Last column shows the convergence rate of $F(x_k) - F^*$ (the optimality gap) and our theorem that proves the claim.
    \begin{table}[H]
        \centering
        {\scriptsize
        \begin{tabular}{|l|l|l|l|l|}
        \hline
            Algorithm & $\mu$ & $\alpha_k$ & $\rho_k$ & Convergence of $F(x_k) - F^*$
        \\ \hline
            R-WAPG in Definition \ref{def:wapg} &
            $\mu \ge 0$ &
            $\alpha_k \in(\mu/L, 1)$ &
            $\rho_k > 0$ &
            \begin{tabular}{l}
                $\mathcal O \left(\prod_{i = 0}^{k-1} \max(1, \rho_i)(1 - \alpha_{i + 1})\right)$
                \\
                (Proposition \ref{prop:wapg-convergence})
            \end{tabular}
        \\ \hline
            Chambolle, Dossal 2015 \cite{chambolle_convergence_2015} &
            $\mu = 0$  &
            $ 0< \alpha_k^{-2} \le \alpha_{k + 1}^{-1} - \alpha_{k + 1}^{-2}$ &
            $\rho_k \ge 1$ &
            \begin{tabular}{l}
                $\mathcal O(\alpha_k^{2})$ \\ (Theorem \ref{thm:r-wapg-on-cham-doss})
            \end{tabular}
        \\ \hline
            V-FISTA Beck (10.7.7) \cite{beck_first-order_2017} &
            $\mu > 0$&
            $\alpha_k = \sqrt{\mu/L}$ &
            $\rho_k = 1$ &
            \begin{tabular}{l}
                $\mathcal O\left((1 - \sqrt{\mu/L})^k\right)$,
                \\
                (Theorem \ref{thm:fixed-momentum-fista}, remark)
            \end{tabular}
        \\ \hline
            R-WAPG in Definition \ref{def:wapg} &
            $\mu > 0$ &
            $\alpha_k = \alpha \in (\mu/L, 1)$ &
            $\rho_k = \rho > 0$ &
            \begin{tabular}{l}
                $\mathcal O \left(\left(1 - \min\left(\mu/(\alpha L), \alpha\right)\right)^{k}\right)$\\
                (Theorem \ref{thm:fixed-momentum-fista})
            \end{tabular}
        \\ \hline
        \end{tabular}
        }
    \end{table}
    The next lemma characterizes momentum sequences in Chambolle and Dossal \cite{chambolle_convergence_2015} using Definition \ref{def:rwapg-seq}.
    \begin{lemma}[R-WAPG sequence as inverted FISTA sequence]\label{lemma:inverted-fista-seq}
        Let R-WAPG sequence $(\rho_k)_{k \ge 0}, (\alpha_k)_{k \ge 0}$ be given by Definition \ref{def:rwapg-seq}.
        If $\mu = 0, \rho_k \ge 1\; \forall k \ge 0$, and $\alpha_0 = 1$, then:
        \begin{enumerate}
            \item $\alpha_k^{-2} \ge \alpha_{k + 1}^{-2} - \alpha_{k + 1}^{-1}\; \forall k \ge 0$.
            \item Let $t_k := \alpha_k^{-1}$, then $0 < t_{k + 1} \le (1/2)\left(1 + \sqrt{1 + 4t_k^2}\right)\;\forall k\ge 0$, hence the name: ``inverted FISTA sequence''.
            \item $\prod_{i = 1}^k\max(1, \rho_{k - 1})(1 - \alpha_k) = \alpha_k^2 \quad (\forall k \ge 1)$.
        \end{enumerate}
    \end{lemma}
    \begin{proof}
        We start proving (i).
        For all $k \ge 1$:
        \begin{align*}
            \alpha_{k + 1}^2
            &= (1 - \alpha_{k + 1})\rho_k\alpha_k^2 + (\mu/L) \alpha_k
            \\
            &= (1 - \alpha_{k + 1})\rho_k\alpha_k^2 &
            (\mu = 0 \text{ assumed } )
            \\
            \implies
            \alpha_{k + 1}^2
            & \ge (1 - \alpha_{k + 1})\alpha_k^2
            &  (\rho_k \ge 1)
            \\
            \iff
            \alpha_k^{-2}
            &\ge
            \alpha_{k + 1}^{-2} - \alpha_{k + 1}^{-1}.
        \end{align*}
        We need to verify the base case for $k = 0$.
        Since $\alpha_0 = 1$, it has:
        \begin{align*}
            \alpha_1^2 &= (1 - \alpha_1)\rho_0\alpha_0^2
            \\
            &\ge (1 - \alpha_1)\alpha_0^2
            = (1 - \alpha_1)
            \\
            \iff
            1 &\ge \alpha_1^{-2} - \alpha_1^{-1}.
        \end{align*}
        Therefore, (i) holds.
        \par
        (i) $\implies$ (ii): Substituting $\alpha_k^{-1} = t_k$ changes (i) into $t_{k + 1}^2 - t_{k + 1} - t_{k}^2 \le 0$ for all $k \ge 0$.
        Solving the equality yields $t_{k + 1} = (1/2)\left(1 \pm \sqrt{1 + 4 t_k^2}\right)$ for all $k \ge 0$.
        Since $\alpha_k \in (0, 1)$, it means $t_k > 0$, hence the valid root is the larger root which gives the upper bound for $t_{k + 1}$ so:
        \begin{align*}
            t_{k + 1} \in \left(
                0, \frac{1}{2}\left(1 + \sqrt{1 + 4t_k^2}\right)
            \right]
        \end{align*}
        To prove (iii), because $\rho_k \ge 1$, we have $\prod_{i = 1}^{k} \max(1, \rho_{k - 1})(1 -
        \alpha_k)= \prod_{i = 1}^{k}\rho_{k - 1}(1 - \alpha_k)$ for all $k \ge 1 $.
        By the definitions of $(\alpha_k)_{k \ge 1}, (\rho_k)_{k \ge 1}$ sequences they have for all $k \ge 1$:
        \begin{align*}
            \alpha_k^2 &= \rho_{k - 1}(1 - \alpha_k)\alpha_{k - 1}^2
            \\
            \iff
            \alpha_k^2/\alpha_{k - 1}^2 &= \rho_{k - 1}(1 - \alpha_k)
            \\
            \implies
            \prod_{i = 1}^{k}\rho_{k - 1}(1 - \alpha_k)
            &=
            \prod_{i = 1}^{k}\alpha_k^2 /\alpha_{k - 1}^2= \alpha_k^2/\alpha_0.
        \end{align*}
        Because $\alpha_0 = 1$, (iii) is justified.
    \end{proof}
    \begin{remark}
        The sequence $(t_k)_{k\geq 1}$ in Lemma~\ref{lemma:inverted-fista-seq} is exactly the same as in Chambolle and Dossal 
        \cite[Theorem 3.1]{chambolle_convergence_2015}.
    \end{remark}
    \begin{theorem}[FISTA first variant Chambolle and Dossal 2015]\label{thm:r-wapg-on-cham-doss}\;\\
        Fix arbitrary $a \ge 2$.
        Define $\forall k \ge 1$ the sequence $(\alpha_k)_{k \ge 0}, (\rho_k)_{k \ge 0}$ by
        \begin{align*}
            \alpha_k &= a/(k + a),
            \\
            \rho_k &= \frac{(k + a)^2}{(k + 1)(k + a + 1)}.
        \end{align*}
        Consider the algorithm given by:
        \begin{tcolorbox}
            Initialize any $y_1 = x_1$.
            \\
            For $k = 1, 2, \ldots$, update:
            \begin{align*}
                & x_{k + 1} := y_k + L^{-1}\mathcal G_L(y_k),
                \\
                & \theta_{k + 1} := \alpha_{k + 1}(\alpha_k^{-1} - 1),
                \\
                & y_{k + 1} := x_{k + 1} + \theta_{k + 1}(x_{k + 1} - x_k).
            \end{align*}
        \end{tcolorbox}
        If $\mu = 0$, then $(\alpha_k)_{k \ge 0}, (\rho_k)_{k \ge 0}$ is a valid pair of R-WAPG sequence from Definition \ref{def:rwapg-seq} and the above algorithm is a valid form of R-WAPG.
        \par
        Assume minimizer $x^*$ exists for function $F$.
        Then algorithm produces $(x_k)_{k \ge 0}$ such that $(F(x_{k}) - F(x^*))_{k\geq 0}$ converges at a rate of $\mathcal O(\alpha_k^2)$.
    \end{theorem}
    \begin{proof}
        Proceeding stepwise, we want to specialize the results of Proposition \ref{prop:wapg-convergence} for the above choice of sequence $(\alpha_k)_{k \ge 0}, (\rho_k)_{k \ge 0}$  by completing the following steps:
        \begin{enumerate}
            \item Show that the $(\alpha_k)_{k \ge 0}, (\rho_k)_{k \ge 0}$ given are valid R-WAPG sequences, so we can use the convergence results.
            \item Show $\rho_k \ge 1$, so Lemma \ref{lemma:inverted-fista-seq} (iii) applies.
            \item Simplify the conclusion of Proposition \ref{prop:wapg-convergence} by the previous two items.
        \end{enumerate}
        To show (i), we verify equality $\alpha_{k + 1}^2 = (1 - \alpha_{k + 1})\rho_k \alpha_k^2$ because it's assumed $\mu = 0$ in this case.
        For all $k \ge 0$, the right side of the equality evaluates to
        \begin{align*}
            \alpha_k^2 \rho_k(1 - \alpha_{k + 1}) &=
            \left(
                \frac{a}{k + a}
            \right)^2 \left(
                \frac{(k + a)^2}{(k + 1)(k + a + 1)}
            \right)
            \left(
                1 - \frac{a}{k + 1 + a}
            \right)
            \\
            &= \left(
                \frac{a}{k + a}
            \right)^2 \left(
                \frac{(k + a)^2}{(k + 1)(k + a + 1)}
            \right)
            \left(
                \frac{k + 1}{k + 1 + a}
            \right)
            \\
            &= \frac{a^2}{(k + a + 1)^2} = \alpha_{k + 1}^2.
        \end{align*}
        Therefore, $(\rho_k)_{k \ge 0}, (\alpha_k)_{k \ge 0}$ is a pair of valid R-WAPG sequence.
        So, it can be used in the R-WAPG algorithm and represent it in R-WAPG Momentum Form in Definition \ref{def:r-wapg-momentum-form}.
        Using $\mu = 0$, it simplifies the momentum term in Definition \ref{def:r-wapg-momentum-form} and gives the rules of updates in the theorem statement.
        \par
        Next, we show (ii).
        Observe that $(\forall k \ge 0)(\forall a \ge 2)$, $\rho_k$ has
        \begin{align*}
            \rho_k &= \frac{(k + a)^2}{(k + 1)(k + a + 1)} = \frac{k^2 + a^2 + 2ak}{k^2 + ak + 2k + a + 1}
            \\
            &= (k^2 + a^2 + 2ak)(k^2 + 2ak - ak + 2k + a + 1)^{-1}
            \\
            &= (k^2 + a^2 + 2ak)(k^2 + 2ak + \underbrace{(2 - a)k + a + 1}_{\le a^2})^{-1} \ge 1.
        \end{align*}
        Hence, (iii) in Lemma \ref{lemma:inverted-fista-seq} applies.
        By Proposition \ref{prop:wapg-convergence}, we have:
        {\small\begin{align*}
            & F(x_{k + 1}) - F(x^*) + \frac{L\alpha_k^2}{2}\Vert v_{k + 1} - x^*\Vert^2
            \\
            & = F(x_{k + 1}) - F(x^*) + \frac{L\alpha_k^2}{2}\Vert x_{k + 1} - x^* + (\alpha_k^{-1} - 1)(x_{k + 1} - x_k)\Vert^2
            \\
            & \le
            \left(
                \prod_{i = 0}^{k - 1} \max(1, \rho_{k})
            \right)
            \left(
                \prod_{i = 1}^{k} \left(1  - \alpha_i\right)
            \right)
            \left(
                F(x_1) - F(x^*) + \frac{L\alpha_0^2}{2}\Vert v_1 - x^*\Vert^2
            \right)
            \\
            &=
            \alpha_k^2
            \left(
                F(x_1) - F(x^*) + \frac{L\alpha_0^2}{2}\Vert v_1 - x^*\Vert^2
            \right)
            & \text{(By Lemma \ref{lemma:inverted-fista-seq}(iii))}
            \\
            &=
            \left(\frac{a}{k + a}\right)^2
            \left(
                F(x_1) - F(x^*) + \frac{L\alpha_0^2}{2}\Vert x_1 - x^*\Vert^2
            \right).
        \end{align*}}
        We can replace $v_1$ to be $x_1$ by Proposition \ref{prop:r-wapg-momentum-repr}.
    \end{proof}
    \begin{remark}
        This algorithm described here is exactly the same algorithm analyzed by Chambolle and Dossal in \cite{chambolle_convergence_2015}.
    \end{remark}
    Following a similar pattern, below we show linear convergences for all $\mu > 0$ strongly convex function with R-WAPG sequence 
    $\alpha_k = \alpha \in (\mu/L, 1)$ that is a fixed constant.

    \begin{lemma}[constant R-WAPG sequence]\label{lemma:constant-rwapg-seq}
        Suppose $(\alpha_k)_{k \ge 0}, (\rho_k)_{k \ge 0}$ are R-WAPG sequences given by Definition \ref{def:rwapg-seq} and assume $L > \mu > 0$.
        Define $q := \mu/L$.
        Then $\forall r \in \left(\sqrt{q},\sqrt{q^{-1}}\right)$, the constant sequence $\alpha_k := r \sqrt{q}$ has the following:
        \begin{enumerate}
            \item For any $r \in \left(\sqrt{q}, \sqrt{q^{-1}}\right)$, the constant sequence $\alpha_k := \alpha \in (q, 1)$ and\\
            $\rho_k := \rho=\left(1-r^{-1}\sqrt{q}\right)\left(1 - r \sqrt{q}\right)^{-1} > 0$, hence it's a pair of valid R-WAPG sequence.
            \item The momentum term $\theta_{k + 1}=\rho_k\alpha_k(1 - \alpha_k)/(\rho_k\alpha_k^2 + \alpha_{k + 1})$ in Definition \ref{def:r-wapg-momentum-form}, denoted by $\theta$ is the constant: $\theta = (1 - r^{-1}\sqrt{q})(1 - r\sqrt{q})(1- q)^{-1}$
            \item When $r = 1$, we have $\theta = (1- \sqrt{q})(1 + \sqrt{q})^{-1}$.
            \item For all $r \in \left(1, \sqrt{q^{-1}}\right)$, we have $\rho > 1$; for all $r \in \left(\sqrt{q}, 1\right]$ we have $\rho \le 1$.
            \item For all $r \in \left(\sqrt{q}, \sqrt{q^{-1}}\right)$, the equality $\max(\rho, 1)(1 - \alpha) = 1 - \min\left(\frac{\mu}{\alpha L}, \alpha\right)$ holds.
        \end{enumerate}
    \end{lemma}
    \begin{proof}
        To see (i), fix any $\sqrt{q} < r < \sqrt{q^{-1}}$, so
        \begin{align*}
            r &\in \left(\sqrt{q}, \sqrt{q^{-1}}\right)
            \iff
            r\sqrt{q} \in
            \left(
                q, 1
            \right).
        \end{align*}
        Therefore, $\alpha_k = \alpha \in (\mu/L, 1)$.
        To see $\rho_k$, by definition it has
        \begin{align*}
            \rho_k &= \frac{\alpha_{k + 1}^2 - q \alpha_{k + 1}}{(1 - \alpha_{k + 1})\alpha_k^2}
            = \frac{\alpha^2 - q \alpha}{(1 - \alpha)\alpha^2}
            \\
            &= \frac{1 - q\alpha^{-1}}{1 - \alpha}
            = \frac{1 - q r^{-1}\sqrt{q^{-1}}}{1 - r \sqrt{q}}
            \\
            &= \frac{1 - r^{-1}\sqrt{q}}{1 - r \sqrt{q}} > 0.
        \end{align*}
        Simple algebra can show (ii).
        To start, we substitute the constant sequence $\alpha, \rho$ for $\alpha_k, \rho_k$ into the definition of $\theta$:
        \begin{align*}
            \theta &= \frac{\rho\alpha(1 - \alpha)}{\rho \alpha^2 + \alpha}
            = (\rho\alpha)\frac{1 - \alpha}{\rho \alpha^2 + \alpha} = \rho \frac{1 - \alpha}{\rho \alpha + 1}
            \\
            &=
            \frac{1 - r^{-1}\sqrt{q}}{1 - r \sqrt{q}}
            (1 - r\sqrt{q})
            \left(
                1 + r\sqrt{q}\frac{1 - r^{-1}\sqrt{q}}{1 - r \sqrt{q}}
            \right)^{-1}
            \\
            &= (1 - r^{-1} \sqrt{q})\left(
                1 + \frac{r \sqrt{q} - q}{1 - r \sqrt{q}}
            \right)^{-1}
            \\
            &= \left(1 - r^{-1} \sqrt{q}\right)\left(
                \frac{1 - q}{1 - r \sqrt{q}}
            \right)^{-1} = \frac{(1 - r^{-1}\sqrt{q})(1 - r \sqrt{q})}{1 - q}.
        \end{align*}
        Now it's the perfect opportunity to show (iii) by substituting $r = 1$ which has
        \begin{align*}
            \theta &= \frac{(1 - \sqrt{q})^2}{1 - q}
            =
            \frac{(1 - \sqrt{q})}{(1 - \sqrt{q})(1 + \sqrt{q})}
            = \frac{1 - \sqrt{q}}{1 + \sqrt{q}}.
        \end{align*}
        Next for (iv), we determine when $\rho$ switches from $> 1$ to $ \le 1$.
        For any $r \in \left(1, \sqrt{q^{-1}}\right)$, we show $\rho > 1$:
        \begin{align*}
            \rho &= \frac{1 - r^{-1}\sqrt{q}}{1 - r \sqrt{q}}
            > \frac{1 - \sqrt{q}}{1 - r \sqrt{q}} > \frac{1 - \sqrt{q}}{1 - \sqrt{q}} = 1
        \end{align*}
        The first inequality comes from $r > 1 \iff -r^{-1} > -1$, the second inequality comes from $r > 1 \iff -r < 1$.
        For any $r \in \left(\sqrt{q}, 1\right]$, we have $\rho \le 1$ by
        \begin{align*}
            \rho &= \frac{1 - r^{-1}\sqrt{q}}{1 - r \sqrt{q}}
            \le \frac{1 - \sqrt{q}}{1 - r \sqrt{q}} \le \frac{1 - \sqrt{q}}{1 - \sqrt{q}} = 1.
        \end{align*}
        The first inequality comes from $-r^{-1} \le - 1$; the second inequality comes from $-r \ge -1$.
        \par
        Finally, to show (v), consider substituting $\rho, \alpha$ in terms of $r$:
        \begin{align*}
            \max(\rho, 1)(1 - \alpha) &=
            \max(\rho(1 - \alpha), (1 - \alpha))
            \\
            &= \max\left(
                \frac{1 - r^{-1}\sqrt{q}}{1 - r \sqrt{q}}(1 - \alpha), (1 - \alpha)
            \right)
            \\
            &=
            \max\left(
                \frac{1 - r^{-1}\sqrt{q}}{1 - r \sqrt{q}}(1 - r\sqrt{q}), (1 - r\sqrt{q})
            \right)
            \\
            &= \max\left(1 - r^{-1}\sqrt{q}, 1 - r \sqrt{q}\right)
            \\
            &= 1 - \min(\sqrt{q}r^{-1}, r\sqrt{q}).
        \end{align*}
        We used the fact that $(1 - \alpha) > 0, \sqrt{q} > 0$ and (i) for substituting $\rho$.
        Now, by definition $r = \alpha/\sqrt{q}$ and $q = \mu/L$, it can be represented by
        \begin{align*}
            1 - \min\left(
                \sqrt{q}\left(\frac{\sqrt{q}}{\alpha}\right), \alpha
            \right) &=
            1 - \min\left(
                q/\alpha, \alpha
            \right) =
            1 - \min\left(\frac{\mu}{\alpha L}, \alpha\right).
        \end{align*}
    \end{proof}

    \begin{theorem}[fixed momentum APG]\label{thm:fixed-momentum-fista}
        Assume $L > \mu > 0$, let a pair of constant R-WAPG sequence: $(\alpha_k)_{k \ge0}, (\rho_k)_{k \ge 0}$ be given by Lemma \ref{lemma:constant-rwapg-seq}.
        Define $q := \mu/L$ and for any fixed $r \in \left(\sqrt{q}, \sqrt{q^{-1}}\right)$, let $\alpha_k := \alpha = r \sqrt{q}$ be the constant R-WAPG sequence.
        Consider the algorithm with a constant momentum specified by the following:
        \begin{tcolorbox}
            Define $\theta = \left(1 - r^{-1}\sqrt{q}\right)(1 - r\sqrt{q})(1 - q)^{-1}$.
            \\
            Initialize $y_1 = x_1$; for $k = 1, 2, \ldots, N$, update:
            \begin{align*}
                &x_{k + 1} = y_k + L^{-1}\mathcal G_L (y_k)
                ,
                \\
                & y_{k + 1} = x_{k + 1} + \theta(x_{k + 1} - x_k).
            \end{align*}
        \end{tcolorbox}
        Then the algorithm generates $(x_k)_{k \ge 1}$ such that $(F(x_{k}) - F(x^*))_{k\geq 1}$ converges at a rate of $\mathcal O\left(\left(1 - \min\left(\mu/(\alpha L), \alpha\right)\right)^k\right)$.
    \end{theorem}
    \begin{proof}
        The constant sequence $(\alpha_k)_{k \ge 0}, (\rho_k)_{k \ge 0}$ is a valid R-WAPG sequence by Lemma \ref{lemma:constant-rwapg-seq}.
        The algorithm presented here is the same as the R-WAPG momentum form as given by Definition \ref{def:r-wapg-momentum-form} because the term $\theta$ from Lemma \ref{lemma:constant-rwapg-seq} matches up with Definition \ref{def:r-wapg-momentum-form} by (ii) in Lemma \ref{lemma:constant-rwapg-seq}.
        Then, Proposition \ref{prop:wapg-convergence} applies, and it follows that:
        \begin{align*}
            & F(x_{k + 1}) - F^* + \frac{L\alpha_k^2}{2}\Vert v_{k + 1} - x^*\Vert^2
            \\
            &=
            F(x_{k + 1}) - F^* + \frac{L\alpha_k^2}{2}\Vert x_{k + 1}  - x^* + (\alpha^{-1} - 1)(x_{k + 1} - x_k)\Vert^2
            \\
            &\le
            \left(
                \prod_{i = 0}^{k - 1} \max(1, \rho_{k})
            \right)
            \left(
                \prod_{i = 1}^{k} \left(1  - \alpha_i\right)
            \right)
            \left(
                F(x_1) - F(x^*) + \frac{L\alpha_0^2}{2}\Vert v_1 - x^*\Vert^2
            \right)
            \\
            &= \left(
                \prod_{i = 1}^{k - 1} \max(1, \rho)(1 - \alpha)
            \right)
            \left(
                F(x_1) - F(x^*) + \frac{L\alpha_0^2}{2}\Vert v_1 - x^*\Vert^2
            \right)
            \\
            &= \left(1 - \min\left(\frac{\mu}{\alpha L}, \alpha\right)\right)^k
            \left(
                F(x_1) - F(x^*) + \frac{L\alpha_0^2}{2}\Vert x_1 - x^*\Vert^2
            \right).
            & \text{(Lemma \ref{lemma:constant-rwapg-seq}(v))}
        \end{align*}
        Going from the second last to the last equality, we used the initialization condition $y_1 = x_1$ specified in Proposition \ref{prop:r-wapg-momentum-repr}, therefore it has $v_1 = x_1$.
    \end{proof}
    \begin{remark}
        When $r = 1$, the algorithm described above is exactly the same as the V-FISTA algorithm specified in (10.7.7) of Beck's book \cite{beck_first-order_2017}.
        The above theorem gives convergence rate for a much relax choice of momentum parameter $\theta$, thus generalizing the convergence of constant step size FISTA to all convex and strongly convex functions and all possible $\alpha_k \in (\mu/L, 1)$.
        In addition, the convergence claim of R-WAPG holds for any $\mu$ smaller than the strong convexity modulus of $f$ that is not zero making the bounds $r \in (\sqrt{q}, \sqrt{q^{-1}})$ looser.
    \end{remark}

\section{The method of Free R-WAPG}\label{sec:free-rwapg}
    In this section, we propose Algorithm \ref{alg:free-rwapg} which estimates the $\mu$ constant as the algorithm executes using the Bregman divergence of $f$ estimated by previous iterates.
    It is called \emph{Free R-WAPG} because it doesn't require any knowledge of $\mu, L$ in prior for the smooth part $f$, making it parameter-free.
    \begin{algorithm}
        \begin{algorithmic}[1]
        {\footnotesize
        \STATE{\textbf{Input: } $f, g, L > \mu \ge 0, x_0 \in \RR^n, N \in \N$}
        \STATE{\textbf{Initialize: }$y_0 := x_0;L_0 := 1; \mu_0 := 1/2; \alpha_0 = 1$;}
        \STATE{\textbf{Compute: } $f(y_k)$; }
        \FOR{$k = 0, 1, 2, \cdots, N$}
            \STATE{\textbf{Compute: }$\nabla f(y_k); x^+:= [I + L_k^{-1}\partial g]^{-1}(y_k - L_k^{-1}\nabla f(y_k))$;}
            \WHILE{$L_k/2\Vert x^+ - y_k\Vert^2 < D_f(x^+, y_k)$}
                \STATE{$L_k:= 2L_k$;}
                \STATE{$x^+ = [I + L_k^{-1}\partial g]^{-1}(y_k - L_k^{-1}\nabla f(y_k))$; }
            \ENDWHILE
            \STATE{$x_{k + 1} := x^+$;}
            \STATE{$\alpha_{k + 1} := (1/2)\left(\mu_k/L_k - \alpha_{k}^2 + \sqrt{(\mu_k/L_k - \alpha_{k}^2)^2 + 4\alpha_{k}^2}\right)$;}
            \STATE{$\theta_{k + 1} := \alpha_k(1 - \alpha_k)/(\alpha_k^2 + \alpha_{k + 1})$;}
            \STATE{$y_{k + 1}:= x_{k + 1} + \theta_{k + 1}(x_{k + 1} - x_k)$; }
            \STATE{\textbf{Compute: } $f(y_{k + 1})$}
            \STATE{$\mu_k := D_f(y_{k + 1}, y_{k})/\Vert y_{k + 1} - y_k\Vert^2 + (1/2)\mu_k$;}
        \ENDFOR
        }
        \end{algorithmic}
        \caption{Free R-WAPG}
        \label{alg:free-rwapg}
    \end{algorithm}
    \par
    In Algorithm \ref{alg:free-rwapg}, lines 5-8 estimate
    upper bound for the Lipschitz constant and find $x^+$, the next iterates produced by proximal gradient descent on previous $y_k$;
    Line 8 updates $x_{k + 1}$ to be $x^+$, a successful iterate identified by the Lipschitz line search routine;
    Line 11 updates the R-WAPG sequence $\alpha_k$ for the iterates $y_{k + 1}$;
    Line 15 updates $\mu$ using the Bregman Divergence of $f$ from iterates $y_{k + 1}, y_k$.
    \par
    Assume $L$ given to initialize the R-WAPG is an upper bound of the Lipschitz smoothness constant of $f$ hence preventing the line search subroutine to be triggered.
    Then the algorithm calls $f(\cdot)$ two times, and $\nabla f(\cdot)$ once per iteration.
    The algorithm computes $\nabla f(y_k)$ once for $x^+$, $f(y_{k + 1})$ once for Bregman Divergence because $f(y_{k})$ is evaluated from the previous iteration, and $f(x^+)$ once for Lipschitz constant line search condition.
    We note that $f(y_0)$ is computed before the start of the for loop.
    And finally, it evaluates proximal of $g$ at $y_k - L^{-1}\nabla f(y_k)$ once.
    \par
    Suppose that initial $L$ given as input to Algorithm \ref{alg:free-rwapg} satisfies $L/2 \ge D_f(x, y)/\Vert x - y\Vert^2\;\forall x, y \in \RR^n$.
    Then, the Lipschitz line search is not triggered and $L_k = L$ for all $k = 1, \ldots, N$.
    Initialize $\mu = (1/2)L$ which it may not be a lower bound on the strong convexity constant for $f$.
    With the additional assumption that the estimated $\mu_k < L$, the sequence $(\alpha_k)_{k \ge 0}$ generated by Algorithm \ref{alg:free-rwapg}.
    Let's verify that $(\alpha_k)_{k \ge 0}$ is an R-WAPG sequence hence the convergence result (Proposition \ref{prop:wapg-convergence}) of R-WAPG applies.
    \par
    Let $\mu_k$ be an estimate produced by the algorithm and assume that it satisfies that $\mu_k/ L \in [0, 1)$ and $\alpha_k \in (\mu/L, 1)$.
    Assuming inductively that $q = \mu_k/L$,  $\alpha_k \in (\mu/L, 1)$, we have
    \begin{align*}
        \alpha_{k + 1}
        &=
        (1/2)\left(
            q - \alpha_k^2 + \sqrt{(q - \alpha_k^2)^2 + 4 \alpha_k^2}
        \right)
        \\
        &=
        (1/2)\left(
            q - \alpha_k^2 + \sqrt{
                (\alpha_k^2 - q + 2)^2 + 4q - 4
            }
        \right)
        \\
        &<
        (1/2)\left(
            q - \alpha_k^2 + \sqrt{
                (\alpha_k^2 - q + 2)^2
            }
        \right) \le 1.
    \end{align*}
    $q \in (0, 1)$ by assumption $\mu < L$ so $4q - 4 < 0$ always.
    It is not hard to see that $\alpha_{k + 1} > \mu/L$ as well.
    Therefore, $\alpha_k \in (0, 1)$.
    To prevent the case $\alpha_k = 1$, we may always force the estimate $\mu = L/2$ whenever $\mu > L/2$.
    \subsection{Numerical experiments}
        This section showcases results for the numerical experiments conducted using FR-WAPG algorithm (Algorithm \ref{alg:free-rwapg}), and compare with other APG algorithms in the literatures: the V-FISTA, and M-FISTA algorithm described in Section (10.7.7, 10.7.6) by Beck \cite{beck_first-order_2017}.
        We implemented in Julia \cite{bezanson_julia_2017}.
        The equivalences highlighted in Proposition \ref{prop:r-wapg-momentum-repr} allows us to compare the sequence of iterates $(x_k)_{k \ge 1}, (y_k)_{k \ge0}$ for R-WAPG, V-FISTA and M-FISTA.
        \par
        We measure the aggregate statistics of the base two logarithms of the normalized optimality gap (NOG), at each iteration $k$ with the same initialization conditions given for all algorithms.
        Given $x_k$, and minimum $F^*$, we define the normalized optimality gap:
        \newcommand{\NOG}{\text{\textbf{NOG}}}
        \begin{align*}
            \delta_k := \log_2\left(
                \NOG_k
            \right)
            \text{ with }
            \NOG_k := \frac{F(x_k) - F^*}{F(x_0) - F^*}.
        \end{align*}
        Since it's not the case that $F^*$ is always known in prior, we will use the minimum of all $F(x_k)$ across all algorithms, all iterations $k$ as the surrogate for $F^*$ when the true value is unavailable.
        \par
        We consider the norm of the gradient mapping $\Vert \mathcal G_L(y_k)\Vert < \epsilon$ as a termination conditions for all test algorithms.
        The $L$ can change during each iteration if it's obtained through the specified Lipschitz line search routine.
        \subsubsection{Simple convex quadratic}
            Consider $\min_{x \in \RR^n} \{F(x):= f(x) + 0\}$ where
            \begin{align*}
                f(x) = (1/2)\langle x, A x\rangle.
            \end{align*}
            $A$ is set to be a positive semi-definite diagonal matrix, so the problem admits the minimizer $x^* = \mathbf 0$ (which is not unique) with minimum being zero.
            We apply Algorithm \ref{alg:free-rwapg}, M-FISTA, and V-FISTA.
            The following parameters are used to set up the numerical experiments:
            \begin{enumerate}
                \item $N$, the dimension of the problem which defines $A \in \RR^{N\times N}$, a diagonal matrix given by $N- 1$ linearly spaced with equal increment on the interval $[\mu, L]$, and an extra number $0$, i.e.,
                \begin{align*}
                    A = \text{diag}
                    \left(0,
                        \mu + \frac{L - \mu}{N - 1},
                        \mu + 2\frac{L - \mu}{N - 1},
                        \ldots,
                        \mu + \frac{(N - 2)(L - \mu)}{N - 1},
                        L
                    \right).
                \end{align*}
                \item The strong convexity and Lipschitz smoothness constant has $0 < \mu < L$. They are given in prior to construct the problem.
                \item $\epsilon > 0$ is the tolerance value, and it's set to be $10^{-10}$.
                \item $x_0 \sim \mathcal N(I, \mathbf 0)$ is a vector, and it's the initial condition for all the algorithm. In this case the initial guess is fixed for all R-WAPG, V-FISTA and M-FISTA, but it's randomly generated by the zero mean standard normal distribution for each element in the vector.
            \end{enumerate}
            The parameter $L=1, \mu=10^{-5}$ are given in prior to produce the diagonal matrix $A$, and we conduct many experiments for $N = 256$ and $N = 1024$.
            For all R-WAPG, M-FISTA and V-FISTA the same random initial guess is used for all test algorithms and 30 experiments are repeated with a different random initial guess each time.
            The maximum, minimum and median values of $\delta_k$ are measured for all algorithms at each iteration and plotted as a ribbon.
            \par
\begin{figure}[H]
                \begin{subfigure}[b]{0.47\textwidth}
                    \centering
                    \includegraphics[width=\textwidth]{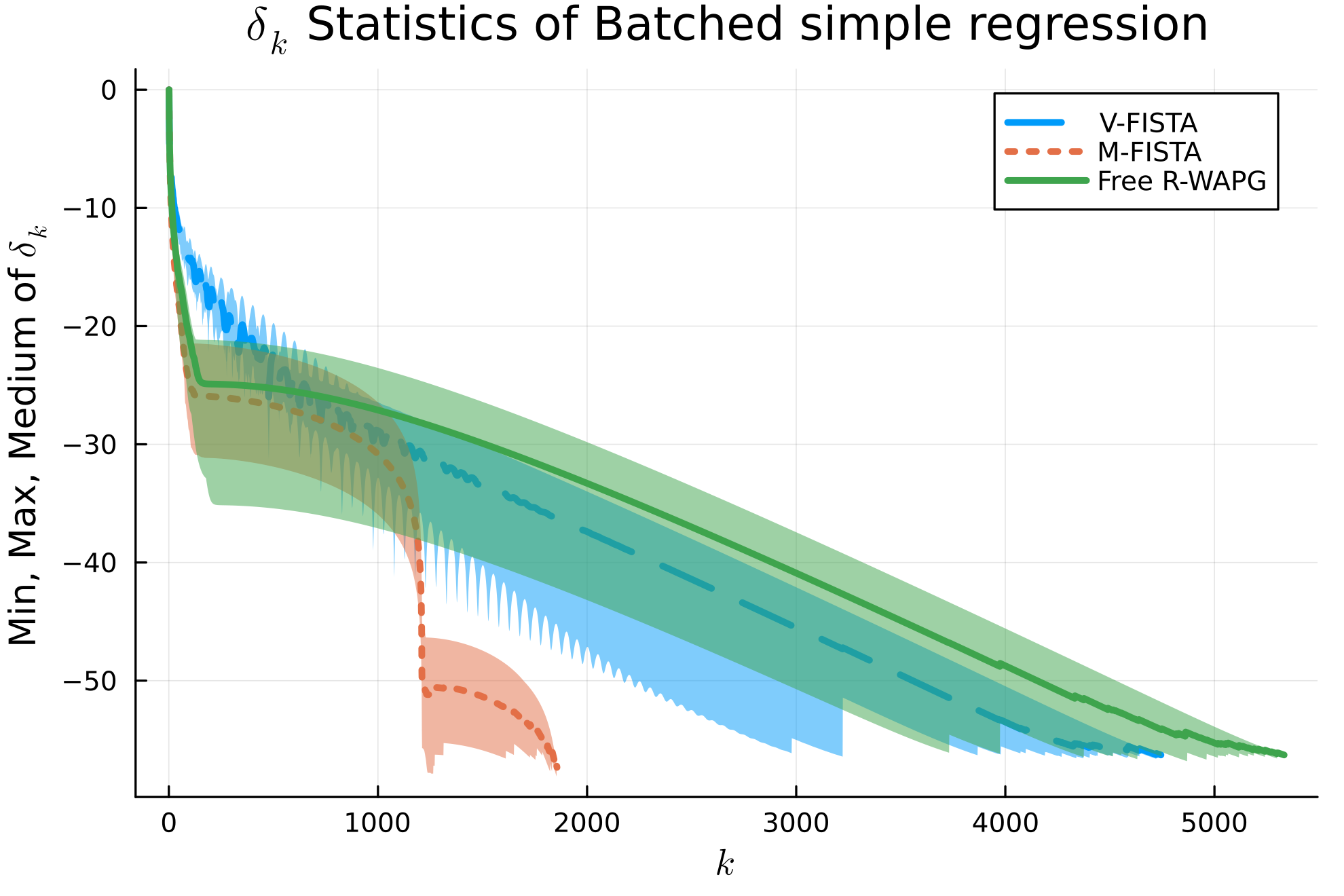}
                    \caption{$N = 256$, simple convex quadratic.}
                \end{subfigure}
                \hfill
                \begin{subfigure}[b]{0.47\textwidth}
                    \centering
                    \includegraphics[width=\textwidth]{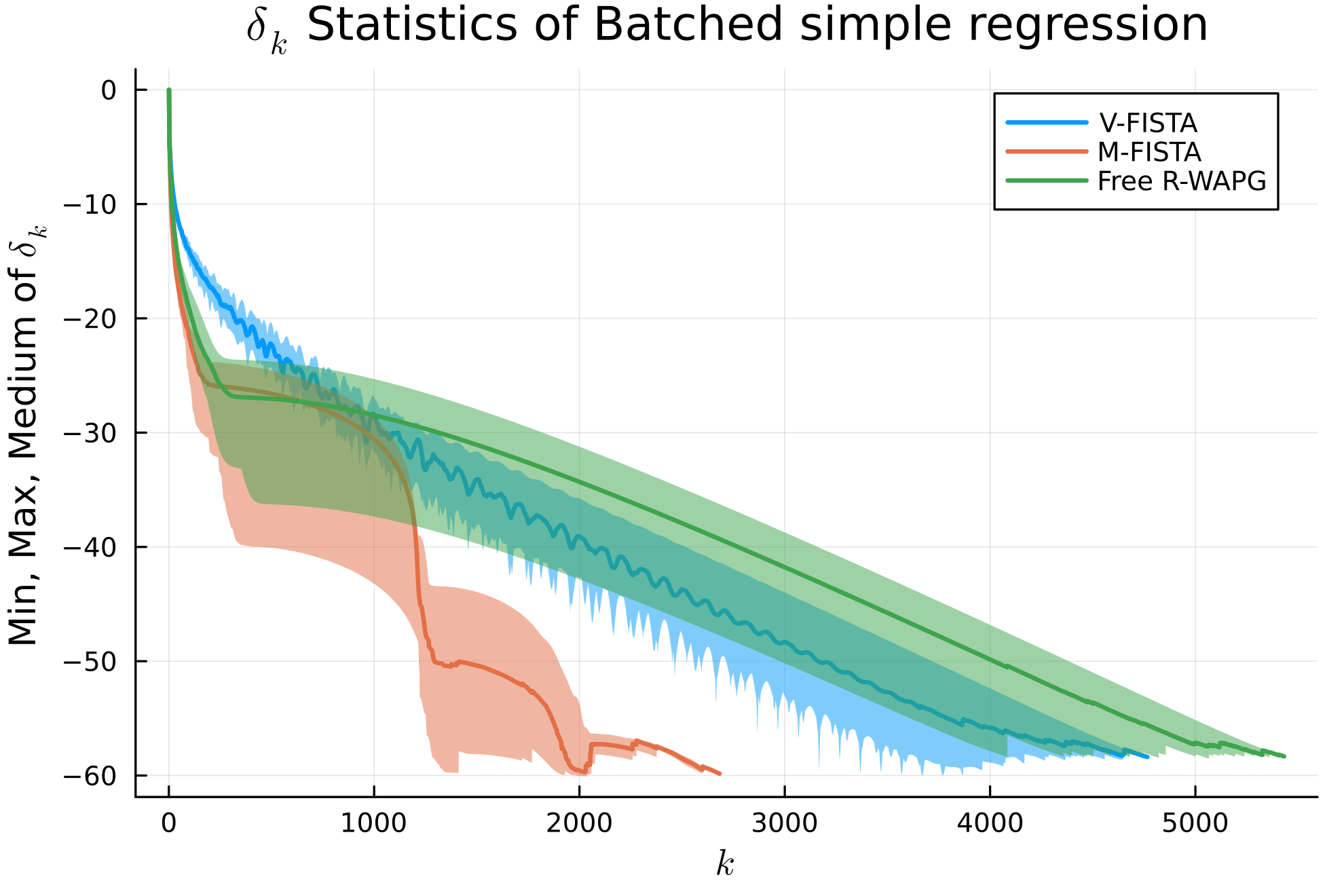}
                    \caption{$N = 1024$, simple convex quadratic. }
                \end{subfigure}
                \caption{
                    Statistics for experiments with simple convex quadratic for V-FISTA, M-FISTA, and R-WAPG.
                }
                \label{fig:simple-quadratic-NOG}
            \end{figure}
\begin{figure}[H]
                \centering
                \begin{subfigure}[b]{0.47\textwidth}
                    \centering
                    \includegraphics[width=\textwidth]{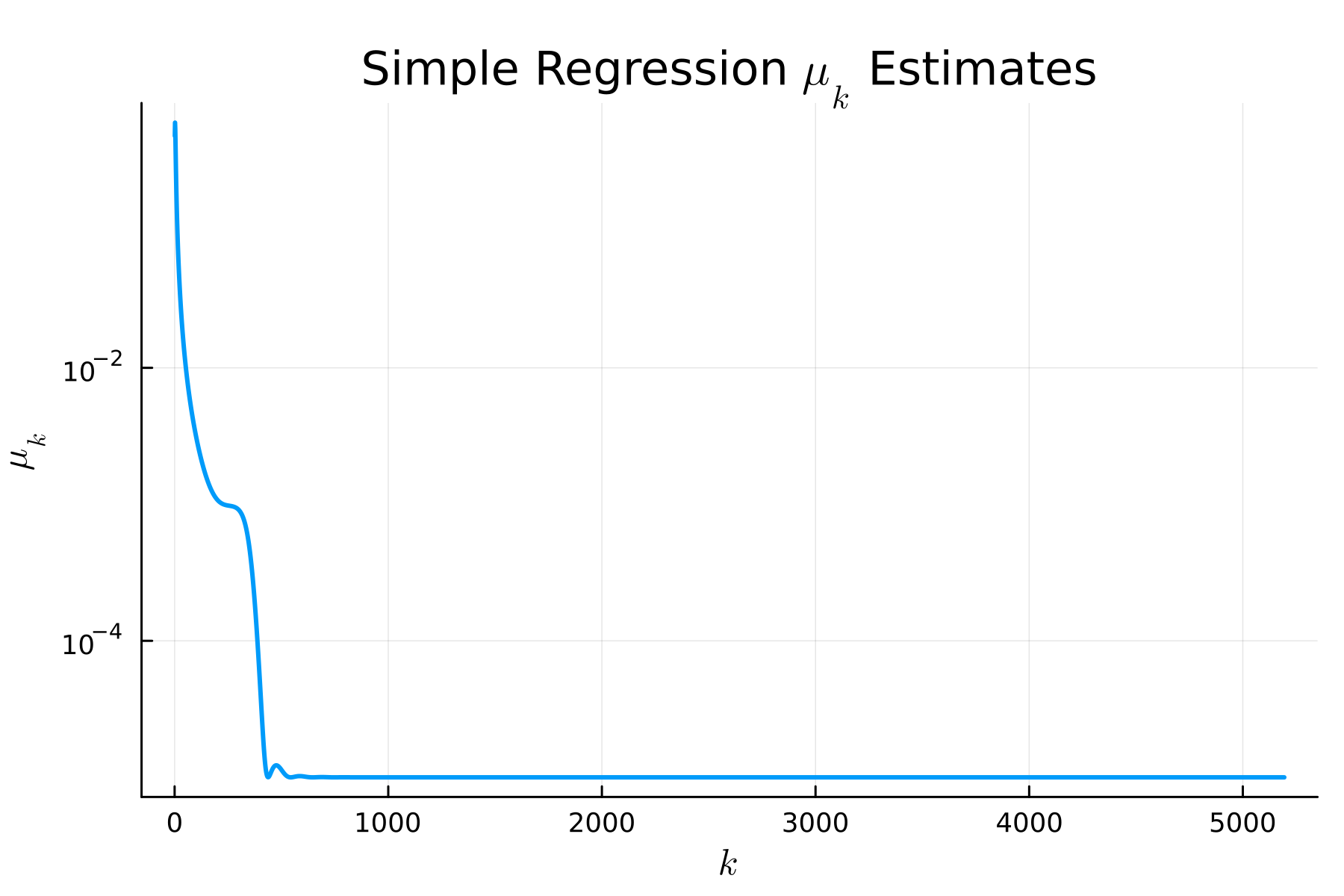}
                \end{subfigure}
                \hfill
                \begin{subfigure}[b]{0.47\textwidth}
                    \centering
                    \includegraphics[width=\textwidth]{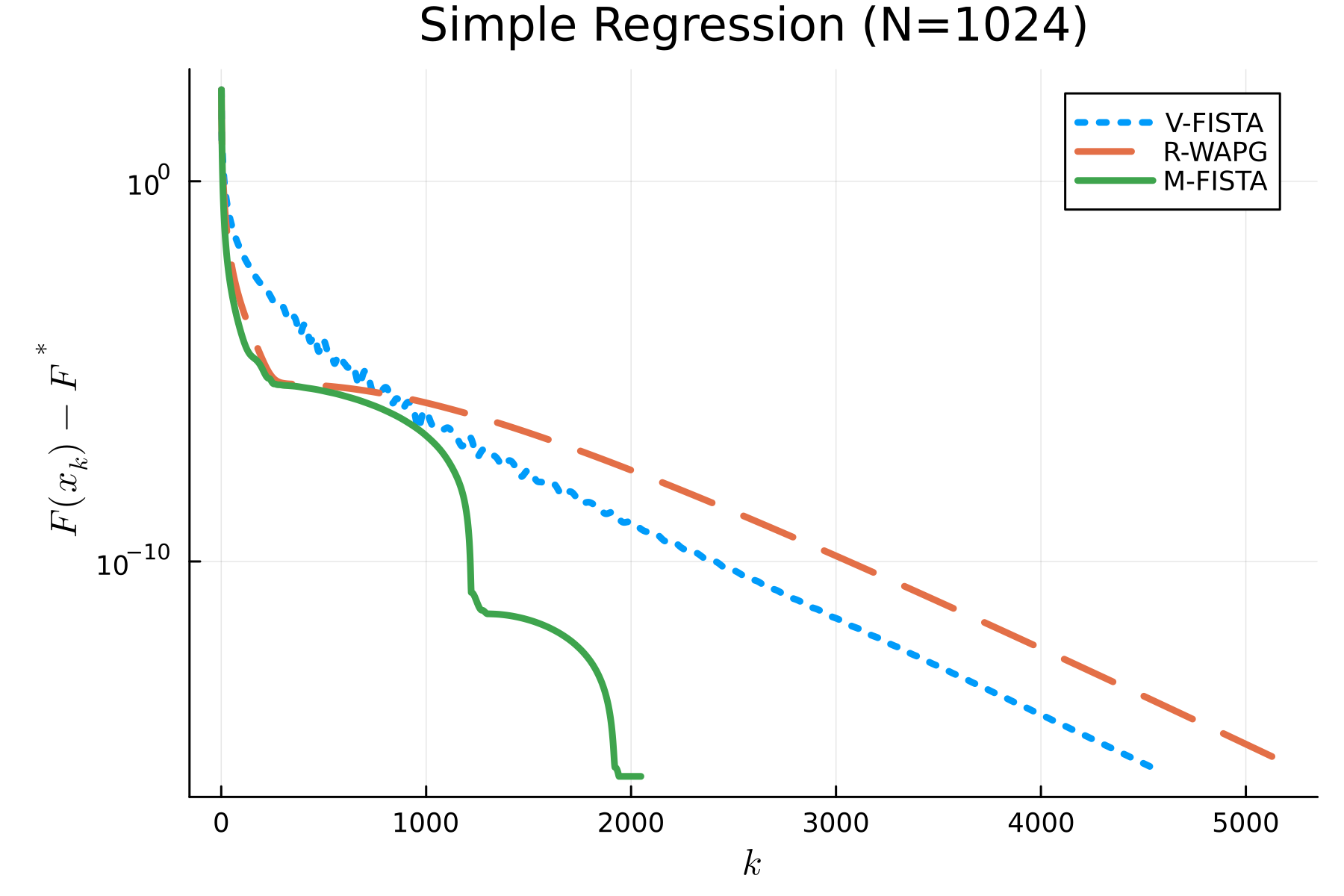}
                \end{subfigure}
                \caption{
                    $N = 1024$, the $\mu$ estimates produced by Algorithm \ref{alg:free-rwapg} (R-WAPG) is recorded.
                }
                \label{fig:simple-quadratic-r-wapg-mu-estimates}
            \end{figure}

            Results are shown in Figure \ref{fig:simple-quadratic-NOG}.
            The solid line in the ribbon is the median value of $\delta_k$ across all experiment, the ribbon gives the maximum, minimum value of $\delta_k$ for each iteration across all experiments.
            This numerical experiment shows that FR-WAPG initially behaves similar to M-FISTA, but it then started to behave like V-FISTA.
            It's numerically verified that the method is not monotone in general, and it just looks monotone in this experiment on this plot.
            Another quantity that maybe interesting other than $\delta_k$ would be the estimated value of $\mu$ during at each iteration $k$.
            One individual experiment is carried out for the R-WAPG algorithm and the value of $\mu$ at each iteration is being recorded as well.
            Figure \ref{fig:simple-quadratic-r-wapg-mu-estimates} shows the estimated $\mu_k$ of FR-WAPG on the left and the convergence rate of the test algorithms on the right for one specific initial condition.
            It shows the values of $\mu_k$ oscillates and converges to the true value of $\mu$.
            Observe that the initial fast descent slows and FR-WAPG diverges from M-FISTA as $\mu_k$ converged.
        \subsubsection{LASSO}
\begin{figure}[H]
                \begin{subfigure}[b]{0.47\textwidth}
                    \centering
                    \includegraphics[width=\textwidth]{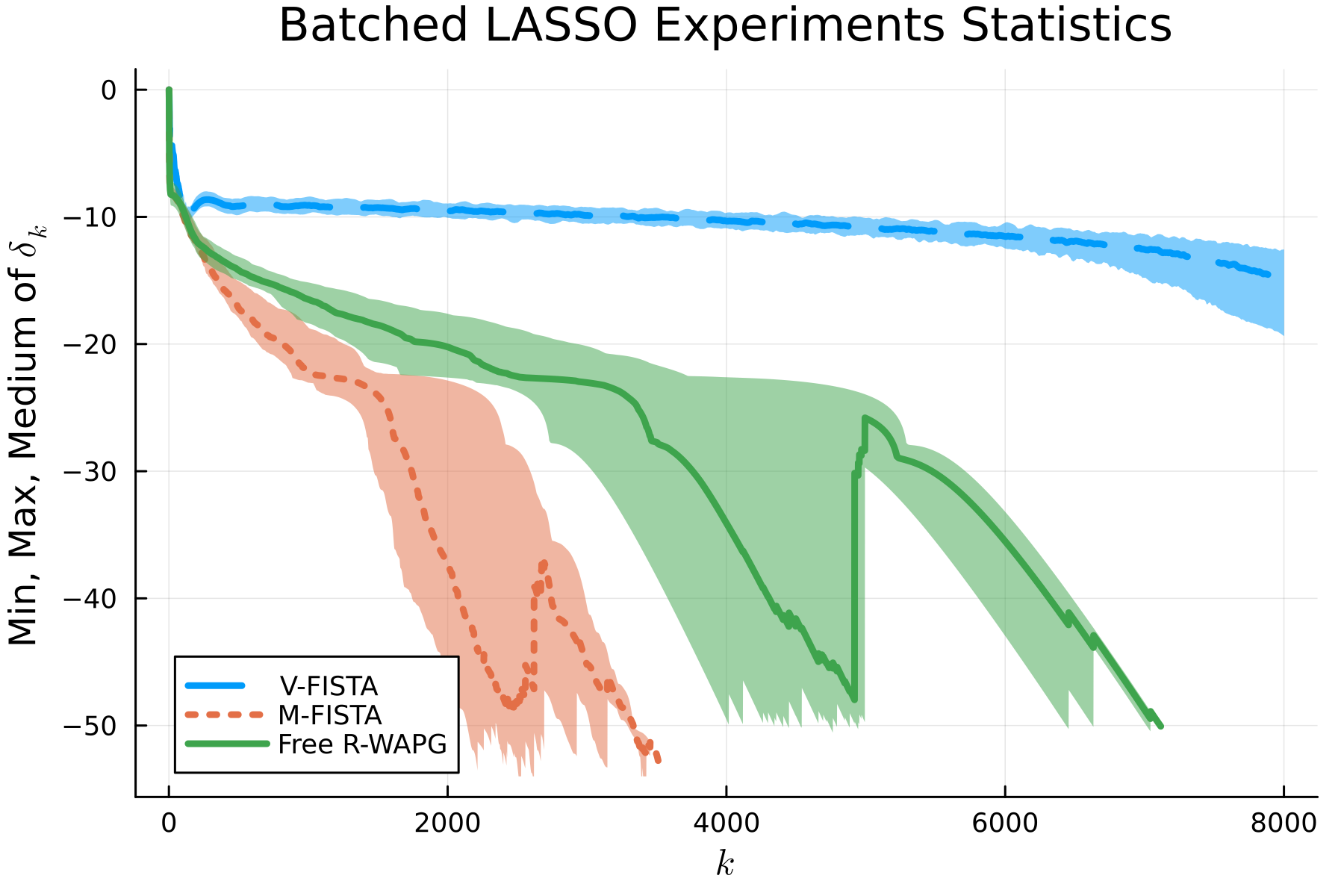}
                    \caption{LASSO experiment with $M = 64, N = 256$. Plots of minimum, maximum, and median $\delta_k$ with estimated $F^*$. }
                \end{subfigure}
                \hfill
                \begin{subfigure}[b]{0.47\textwidth}
                    \centering
                    \includegraphics[width=\textwidth]{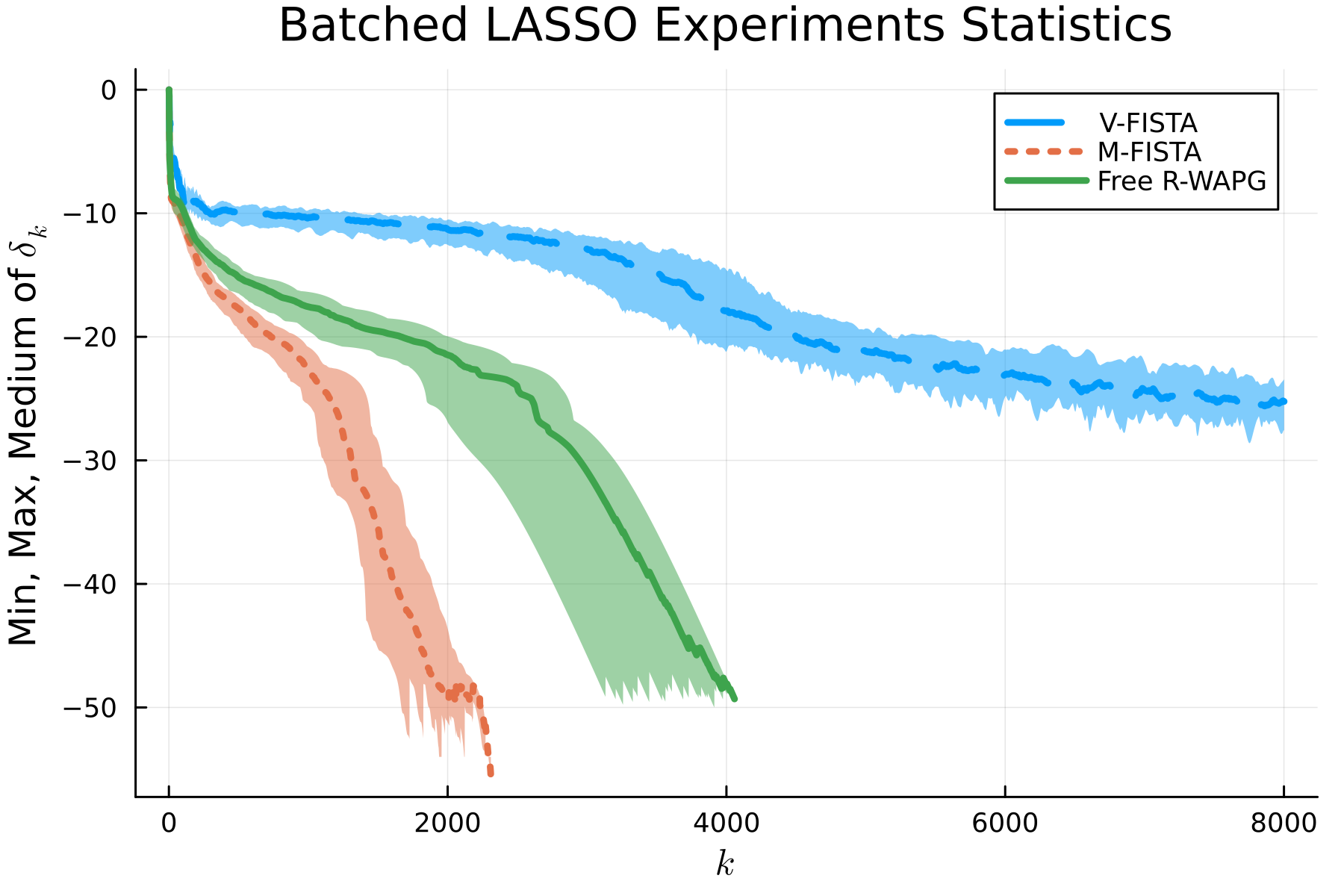}
                    \caption{LASSO experiment with $M = 64, N = 128$. Plots of minimum, maximum, and median $\delta_k$ with estimated $F^*$. }
                \end{subfigure}
                \caption{LASSO experiments statistics for test algorithms. }
                \label{fig:batched-lasso}
            \end{figure}

\begin{figure}[H]
                \begin{subfigure}[b]{0.47\textwidth}
                    \includegraphics[width=\textwidth]{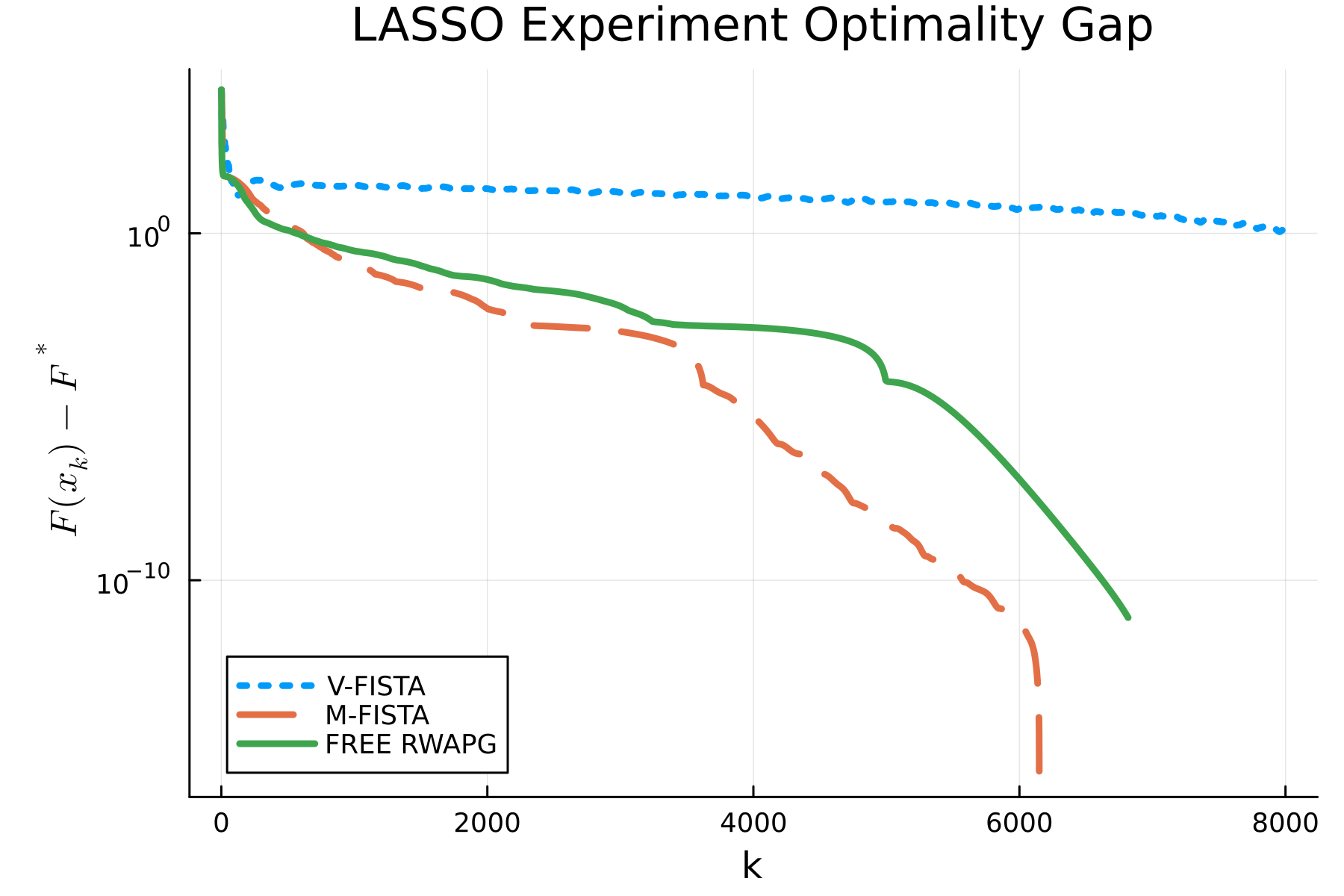}
                    \caption
                    {A single run of LASSO experiment displaying $F(x_k) - F^*$ for several test algorithms.
                    }
                \end{subfigure}
                \hfill
                \begin{subfigure}[b]{0.47\textwidth}
                    \includegraphics[width=\textwidth]{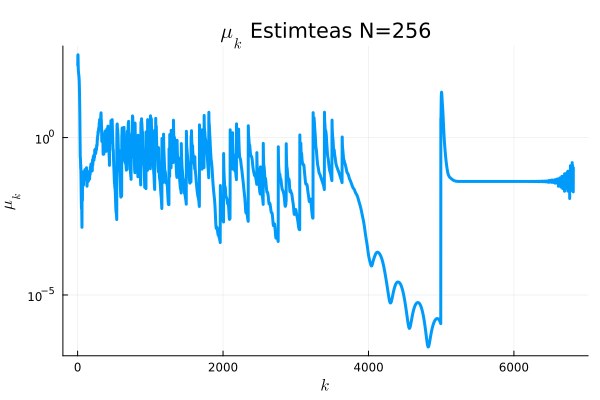}
                    \caption{The $\mu_k$ estimated by FR-WAPG for one LASSO experiment. }
                \end{subfigure}
                \caption{A single LASSO experiment results, with $M = 64, N=256$. f}
                \label{fig:single-lass-mu-estimates}
            \end{figure}

            This section present results of numerical experiment for solving the LASSO (least absolute shrinkage and selection operator) problem proposed by Tibshirani \cite{tibshirani_regression_1996}.
            The problem of LASSO has smooth and nonsmooth additive parts given by:
            \begin{align*}
                \min_{x \in \RR^n}
                \left\lbrace
                    \frac{1}{2}\Vert Ax - b\Vert^2 + \lambda\Vert x\Vert_1
                \right\rbrace.
            \end{align*}
            The smooth part is $f(x) :=\frac{1}{2}\Vert Ax - b\Vert^2$ and the nonsmooth is $g(x) := \lambda\Vert x\Vert_1$.
            The objective function is coersive and the exact minimum, or minimizers are unknown.
            We perform numerical experiments using V-FISTA, M-FISTA and FR-WAPG on this problem.
            The setup of our parameters now follow.
            \begin{enumerate}
                \item $M, N$ are constants. They define matrix $A \in \RR^{M\times N}$ which has entries of i.i.d random variable taken from a standard normal distribution.
                \item $L, \mu$, are Lipschitz constant and the strong convexity constant for the smooth part of the objective which are not known prior. Consequently, they are estimated by $A$ by $\mu = 1/\Vert (A^TA)^{-1}\Vert$ and $L = \Vert A^TA\Vert$ prior to experiment to assist the V-FISTA algorithm.
                \item $x^+ = [1,\; -1,\; 1, \; \cdots ]^T \in \RR^N$, it's a vector with alternating $1, -1$ in it.
                \item Given $x^+$, it has $b = Ax^+ \in \RR^M$.
                \item Given $A$, estimations for $L,\mu$ are given by $L = \Vert A^TA\Vert$, $\mu = \Vert (A^TA)^{-1}\Vert^{-1}$.
                \item $x_0\in \RR^N$ denotes the initial guess consist of i.i.d random variable realized from the standard normal distribution.
                \item $\epsilon = 10^{-6} > 0$ sets the tolerance for all test algorithm on $\Vert \mathcal G_L(y_k)\Vert$.
            \end{enumerate}
            Experiments were conducted using V-FISTA, M-FISTA and FR-WAPG with $(M, N) = (64, 256)$ and $(M, N) = (64, 128)$ respectively.
            A random matrix $A$ is realized and is fixed for all test algorithms and all repetitions.
            The same experiment are repeated 30 times on all algorithms with a fixed random initial condition $x_k$ realized from the distribution.
            The aggregate statistics of $\delta_k$ are collected for all repetitions, and then grouped by the respective algorithm.
            Figure \ref{fig:batched-lasso} shows the results of the two experiments.
            The bump on the curve is due to a subset of test instances of the 30 repetition where the algorithms take larger number of iterations to terminate.
            The experiment results shows faster convergence rate of FR-WAPG relative to V-FISTA, and it also shows that it behaves similarly to M-FISTA.
            In addition, an experiment were conducted to show the estimate $\mu_k$ and optimality gap $F(x_k) - F^*$ in Figure \ref{fig:single-lass-mu-estimates}.
            For this specific experiment showed in the figure, the estimated value of $\mu, L$ which we feed into V-FISTA are $\mu = 7.432363627613958\times 10^{-18}$ and $L = 2321.737206983643$.
            An important observation is that $(\mu_k)_{k\geq 1}$ didn't converge to the true value in this case, however it didn't affect the convergence of $(F(x_k) - F^*)_{k\geq 1}$.

\section{Discussions}\label{s:discussion}
    This paper provides upper bounds on the optimality gap $F(x_k) - F^*$ for smooth plus nonsmooth composite optimization problems with convex objectives under much weaker assumptions for the momentum parameter.
    The proposed R-WAPG unifies the convergence proof of major Euclidean variants of FISTA with strong convexity parameter $\mu \ge 0$.
    To the best of our knowledge, there is no framework had been identified yet in the literature that proves the convergences of several variants (with and without strong convexity parameter) with momentum sequences that doesn't follow Nesterov's rule.
    \par
    In addition, we proposed explorative numerical experiments for Free R-WAPG showcasing a parameter free algorithm that doesn't use the idea of restarting to achieve competitive convergence rate.
    Unfortunately a precise description of the convergence rate of Free R-WAPG remains a mystery.
    It will be interesting to explore this line of research.
    Our convergence framework gives an upper bound using the sequence $(\alpha_k)_{k \ge 0}$ in $(\mu/L, 1)$, but it requires knowing what the sequence would be in advanced given the initial conditions and the objective function.
    \par
    From a practical perspective, Free R-WAPG (Algorithm \ref{alg:free-rwapg}) shows competitive results of convergence for strongly convex functions.
    The algorithm behaves in a way having characteristics of both V-FISTA and M-FISTA which is fascinating.
    In addition, of being parameter free, the algorithm works well with auto-differentiation framework because it requires the function value at the same point where the gradient is evaluated at $y_k$.
    We do note that in the case when function is not strongly convex, the function value converges without $\Vert \mathcal G_L(y_k)\Vert$ converging.
    This will be explored more in future work.

\section*{Acknowledgements}
The research of HL and XW was partially supported by
the NSERC Discovery Grant of Canada.

\bibliographystyle{siam}
\bibliography{references/refs.bib}

\begin{thebibliography}{10}

\bibitem{ahn_understanding_2022}
{\sc K.~Ahn and S.~Sra}, {\em Understanding {Nesterov}'s acceleration via proximal point method}, in Symposium on {Simplicity} in {Algorithms}, SIAM, June 2022, pp.~117--130.

\bibitem{apidopoulos_convergence_2018}
{\sc V.~Apidopoulos, J.-F. Aujol, and C.~H. Dossal}, {\em Convergence rate of inertial {Forward}-{Backward} algorithm beyond {Nesterov}'s rule}, Mathematical Programming, Series A,  (2018), pp.~1--20.

\bibitem{attouch_first-order_2022}
{\sc H.~Attouch, Z.~Chbani, J.~Fadili, and H.~Riahi}, {\em First-order optimization algorithms via inertial systems with {Hessian} driven damping}, Mathematical Programming, 193 (2022), pp.~113--155.

\bibitem{aujol_parameter-free_2024}
{\sc J.-F. Aujol, L.~Calatroni, C.~Dossal, H.~Labarrière, and A.~Rondepierre}, {\em Parameter-{Free} {FISTA} by adaptive restart and backtracking}, SIAM Journal on Optimization, 34 (2024), pp.~3259--3285.

\bibitem{beck_first-order_2017}
{\sc A.~Beck}, {\em First-order {Methods} in {Optimization}}, {MOS}-{SIAM} {Series} in {Optimization}, SIAM, 2017.

\bibitem{beck_fast_2009}
{\sc A.~Beck and M.~Teboulle}, {\em Fast gradient-based algorithms for constrained total variation image denoising and deblurring problems}, IEEE Transactions on Image Processing, 18 (2009), pp.~2419--2434.

\bibitem{bezanson_julia_2017}
{\sc J.~Bezanson, A.~Edelman, S.~Karpinski, and V.~B. Shah}, {\em Julia: {A} fresh approach to numerical computing}, SIAM Review, 59 (2017), pp.~65--98.

\bibitem{chambolle_convergence_2015}
{\sc A.~Chambolle and C.~Dossal}, {\em On the convergence of the iterates of the ``{Fast} iterative shrinkage/thresholding algorithm"}, Journal of Optimization Theory and Applications, 166 (2015), pp.~968--982.

\bibitem{lee_geometric_2021}
{\sc J.~Lee, C.~Park, and E.~Ryu}, {\em A {Geometric} structure of acceleration and its role in making gradients small fast}, in Advances in {Neural} {Information} {Processing} {Systems}, vol.~34, 2021, pp.~11999--12012.

\bibitem{maulen_speed_2023}
{\sc J.~J. Maulén and J.~Peypouquet}, {\em A speed restart scheme for a dynamics with hessian-driven damping}, Journal of Optimization Theory and Applications, 199 (2023), pp.~831--855.

\bibitem{necoara_linear_2019}
{\sc I.~Necoara, Y.~Nesterov, and F.~Glineur}, {\em Linear convergence of first order methods for non-strongly convex optimization}, Mathematical Programming, 175 (2019), pp.~69--107.

\bibitem{nesterov_method_1983}
{\sc Y.~Nesterov}, {\em A method for solving the convex programming problem with convergence rate {O}(1/k{\textasciicircum}2)}, Proceedings of the USSR Academy of Sciences,  (1983).

\bibitem{nesterov_lectures_2018}
{\sc Y.~Nesterov}, {\em Lectures on {Convex} {Optimization}}, vol.~137 of Springer {Optimization} and {Its} {Applications}, Springer International Publishing, 2018.

\bibitem{ryu_large-scale_2022}
{\sc E.~K. Ryu and W.~Yin}, {\em Large-scale {Convex} {Optimization}: {Algorithms} \& {Analyses} via {Monotone} {Operators}}, Cambridge University Press, 2022.

\bibitem{su_differential_2016}
{\sc W.~Su, S.~Boyd, and E.~J. Candes}, {\em A differential equation for modeling nesterov's accelerated gradient method: {Theory} and {Insights}}, Journal of Machine Learning Research, 17 (2016), pp.~1--43.

\bibitem{tibshirani_regression_1996}
{\sc R.~Tibshirani}, {\em Regression shrinkage and selection via the {Lasso}}, Journal of the Royal Statistical Society. Series B (Methodological), 58 (1996), pp.~267--288.

\end{thebibliography}

\appendix

\end{document}